\documentclass{amsart}
\usepackage{times, pdfsync}

\usepackage{xypic}
\input xy
\xyoption{all}
\usepackage{epsfig}
\usepackage{amsthm}
\usepackage{amssymb}
\usepackage{amsmath}
\usepackage{amscd}
\usepackage{color}

\newcommand{\bg}{\begin{equation}}
\newcommand{\ed}{\end{equation}}
\newcommand{\bga}{\begin{eqnarray}}
\newcommand{\eda}{\end{eqnarray}}
\newcommand{\pf}{\textbf{Proof:\ }}
\newcommand{\lpf}{\textbf{Proof of Lemma:\ }}

\def\cbdu{\par{\raggedleft$\Box$\par}}

\newtheorem {Theorem}  {Theorem}

\numberwithin{Theorem}{section}

\newtheorem {Lemma}[Theorem]  {Lemma}

\theoremstyle{definition}
\newtheorem{Definition}[Theorem]{Definition}
\theoremstyle{remark}
\newtheorem{Remark}[Theorem]{\bf Remark}

\newtheorem {Corollary}[Theorem]{\bf Corollary}
\newtheorem{Claim}[Theorem]{Claim}
\newtheorem{Assumption}[Theorem]{Assumption}
\expandafter\chardef\csname pre amssym.def
at\endcsname=\the\catcode`\@ \catcode`\@=11
\def\undefine#1{\let#1\undefined}
\def\newsymbol#1#2#3#4#5{\let\next@\relax
 \ifnum#2=\@ne\let\next@\msafam@\else
 \ifnum#2=\tw@\let\next@\msbfam@\fi\fi
 \mathchardef#1="#3\next@#4#5}
\def\mathhexbox@#1#2#3{\relax
 \ifmmode\mathpalette{}{\m@th\mathchar"#1#2#3}%
 \else\leavevmode\hbox{$\m@th\mathchar"#1#2#3$}\fi}
\def\hexnumber@#1{\ifcase#1 0\or 1\or 2\or 3\or 4\or 5\or 6\or 7\or 8\or
 9\or A\or B\or C\or D\or E\or F\fi}

\font\teneufm=eufm10 \font\seveneufm=eufm7 \font\fiveeufm=eufm5
\newfam\eufmfam
\textfont\eufmfam=\teneufm \scriptfont\eufmfam=\seveneufm
\scriptscriptfont\eufmfam=\fiveeufm

\catcode`\@=\csname pre amssym.def at\endcsname

\newcounter{remark}
\newcommand{\e}{\epsilon}

\newcommand{\R}{\mathbf{R}}
\renewcommand{\div}{\mbox{div}}
\newcommand{\Aa}{{\mathcal A}}
\newcommand{\Bb}{{\mathcal B}}
\newcommand{\Cc}{{\mathcal C}}

\newcommand{\Dd}{{\mathcal D}}

\newcommand{\Ll}{{\mathcal L}}

\newcommand{\Hh}{{\mathcal H}}

\newcommand{\Mm}{{\mathcal M}}

\def  \R   {{\mathbb R}}

\def  \12  {{\frac{1}{2}}}

\def\build#1_#2^#3{\mathrel{\mathop{\kern 0pt#1}\limits_{#2}^{#3}}}

\begin{document}
%\currannalsline{0}{2006}

\title[Full Liquid Crystal system]{Existence of regular solutions to the full Liquid Crystal System}

\author[Mimi Dai]{ Mimi Dai}
\address{Department of Applied Mathematics, University of Colorado Boulder, Boulder, CO, 80309,USA}
\email{mimi.dai@colorado.edu}

%\shorttitle{Steady-state Navier-Stokes in $\mathbb{R}^3$}
%for running head if title too long.  Comment out if not needed

%%%%%%%%%%The body of the paper follows.
%\thanks{The work of Mimi Dai was partially supported by NSF Grant
%DMS-0700535 and DMS-0900909;}
%%%use \Proof instead of \begin{proof}
%%%% use \Endproof instead of \end{proof}
%%%% use \references {999} instead of \begin{thebibliography}{99}
%%%%used \Endrefs instead of \end{thebibliography}

\begin{abstract}
We study the general Ericksen-Leslie system with non-constant density, which describes the flow of nematic liquid crystal. In particular the model investigated here is associated with Parodi's relation. 
%We first establish the existence of global weak solutions. Next, with more regular initial data,  
We prove that: in two dimension, the solutions are globally regular with general data; in three dimension, the solutions are globally regular with small initial data, or for short time with large data.  Moreover, a weak-strong type of uniqueness result is obtained.  
\medskip 
 
\hspace*{0.4cm}{\it Classification code: 76D03, 3Q35.}
\end{abstract}

\maketitle

%888888888888888888888888888888888888888888888888888888888888888888888888
\section{Introduction}

The flows of nematic liquid crystals can be treated as slow moving particles where the fluid velocity and the alignment of the particles influence each other. The hydrodynamic theory of liquid crystals  was established by Ericksen \cite{Er0, Er1} and Leslie \cite{Le0, Le1} in the 1960's.
As Leslie points out in his 1968 paper : ``liquid crystals are states of matter which are capable of flow, and in which the molecular arrangements  give rise to a preferred direction".  The full Ericksen-Leslie system consists of the following equations (cf. \cite{Er0, Le0, Le1, LL1}):

\begin{equation}\label{LCD}
\begin{split}
\rho_t+&\nabla\cdot (\rho u) =0,\\
\rho \dot u =&\rho F+\nabla\cdot\hat\sigma,\\
\rho_1\dot\omega=&\rho_1 G+\hat g+\nabla\cdot\pi,
\end{split}
\end{equation}
in $\Omega\times(0, T)$, where $\Omega$ is a domain in $R^n$ with $n=2,3$. The three equations in system (\ref{LCD}) describe the conservation of mass, linear momentum and angular momentum, respectively. The anisotropic feature of liquid crystal materials is exhibited in the third equation and the nonlinear coupling is represented in second equation.
%Regularity will be obtained in the following cases,
%\begin{enumerate}
%\item $\Omega=\R^2$ or $\Omega\subset\R^2$ a smooth bounded domain,
%\item $\Omega=\R^3$ or $\Omega\subset\R^3$ a smooth bounded domain, for short time,
%\item $\Omega=\R^3$ or $\Omega\subset\R^3$ a smooth bounded domain, for small data.
%\end{enumerate}
In the above equations,
$\rho: \Omega\times[0,T]\to\mathbb{R}$ is the fluid density, 
$u: \Omega\times[0,T]\to\mathbb{R}^n$ is the fluid velocity,
$d: \Omega\times[0,T]\to\mathbb{R}^n$ is the director field representing the alignment of the molecules,  $\rho_1$ is the inertial constant, $\hat g$ is the intrinsic force associated with $d$, $\pi$ is the director stress, $F$ and $G$ are external body force and external director body force, respectively. 
In this paper we consider the incompressible flow with $\nabla\cdot u=0$. 
The superposed dot denotes the material derivative $\partial_t+u\cdot\nabla$. The notations
\begin{align}\notag
&A=\frac{1}{2}(\nabla u+\nabla^T u), \ \ \ \ \Omega=\frac{1}{2}(\nabla u-\nabla^T u),\\
&\omega=\dot d=d_t+(u\cdot\nabla)d, \ \ \ \ N=\omega-\Omega d,\notag
\end{align}
represent the rate of the strain tensor, the skew-symmetric part of the strain rate, the material derivative of $d$ and the rigid rotation part of director changing rate by fluid vorticity, respectively.
%$d: \Omega\times[0,T]\to\mathbb{R}^n$ is the director field representing the alignment of the molecules, with $n=2,3$. The force term $\nabla d\otimes\nabla d$ in the equation of the conservation of momentum denotes the $3\times 3$ matrix whose $ij$-th entry is given by $``\nabla_i d\cdot\nabla_j d"$ for $1\leq i,j\leq 3$. This force $\nabla d\otimes\nabla d$ is the stress tensor of the energy about the director field $d$, where the energy is given by:
%$$
%\frac 12 \int_{\Omega} |\nabla d|^2 dx + \int_{\Omega}F(d)dx
%$$
%and
%$$
%F(d)=\frac{1}{4\eta^2}(|d|^2-1)^2, \quad f(d) = \nabla F(d) = \frac{1}{\eta^2}(|d|^2-1)d.
%$$
%Hence
%\begin{equation}
%f(d)=\nabla_dF(d)=\frac{1}{\eta^2}(|d|^2-1)d.\notag
%\end{equation}
%In fact $F(d)$ is the penalty term of the Ginzburg-Landau approximation of the original free energy of the director field with unit length.\\

We have the following constitutive relations for $\hat\sigma, \pi$ and $\hat g$ in (\ref{LCD}):
\begin{equation}\label{re:sigma}
\begin{split}
%\label{re:sigma}
\hat\sigma_{ij}=&-P\delta_{ij}-\frac{\partial(\rho W)}{\partial{d_{k,i}}}d_{k,j}+\sigma_{ij},\\
%\label{re:pi}
\pi_{ij}=&\beta_id_j+\frac{\partial(\rho W)}{\partial{d_{j,i}}},\\
%\label{re:g}
\hat g_{i}=&\gamma d_i-\beta_jd_{i,j}-\frac{\partial(\rho W)}{\partial{d_{i}}}+g_{i}.
\end{split}
\end{equation}
Here the scalar function $P$ represents the pressure. The vector $\beta=(\beta_1,\beta_2,\beta_3)^T$ and the scalar function $\gamma$ are Lagrangian multipliers for the constraint $|d|=1$. The term $\rho W$ denotes the Oseen-Frank energy functional for the equilibrium configuration of a unit director field. For simplicity, we consider the relaxation form of the elastic energy associated with $d$:
\bg\label{eq:W}
\rho W=\frac{1}{2}|\nabla d|^2+\frac{1}{4\eta^2}(|d|^2-1)^2
\ed
with constant $\eta>0$. And
\begin{equation}\label{eq:g}
\begin{split}
g_i=&\lambda_1N_i+\lambda_2 d_jA_{ji},\\
\sigma_{ij}=&\mu_1d_kA_{kp}d_pd_id_j+\mu_2N_id_j+\mu_3d_iN_j\\
&+\mu_4A_{ij}+\mu_5A_{ik}d_kd_j+\mu d_iA_{jk}d_k,\\
N_i=&\omega_i+\Omega_{ki}d_k. 
\end{split}
\end{equation}
In the above equations, the constants $\lambda_1, \lambda_2$ represent the molecular shape and $\mu_1, \ldots, \mu_6$ are Leslie coefficients which regard certain local correlations in the fluid (cf. \cite{GP}). They satisfy (cf. \cite{Le0})
\begin{align}
&\lambda_1=\mu_2-\mu_3, \ \ \ \lambda_2=\mu_5-\mu_6, \label{re:con}\\
&\mu_5-\mu_6=-(\mu_2+\mu_3). \label{eq:Par}
\end{align}
The relations in (\ref{re:con}) arise from the second law of thermodynamics. The relation (\ref{eq:Par}) is called Parodi's condition, which is derived from the Onsager's reciprocal relation.\\

To further simplify the model, we take $\rho_1=0$, $\beta=0$, $\gamma=0$, and $F=G=0$.
Thus, the incompressible Ericksen-Leslie system (\ref{LCD}) is reformulated as
\begin{align}
&\rho_t+(u\cdot\nabla) \rho=0, \label{eqL1}\\
&\rho u_t+\rho (u\cdot\nabla) u+\nabla P=-\nabla\cdot(\nabla d\otimes\nabla d)+\nabla\cdot\sigma, \label{eqL2}\\
&d_t+(u\cdot\nabla)d-\Omega d+\frac{\lambda_2}{\lambda_1}Ad=-\frac{1}{\lambda_1}(\Delta d-f(d)), \label{eqL3}\\
&\nabla\cdot u=0, \label{eqL4}
\end{align}
with
$$
\rho\mathcal F(d)=\frac{1}{4\eta^2}(|d|^2-1)^2, \quad f(d) = \nabla_d (\rho \mathcal F(d)) = \frac{1}{\eta^2}(|d|^2-1)d.
$$
The force term $\nabla d\otimes\nabla d$ in the equation of the conservation of momentum denotes the $3\times 3$ matrix whose $ij$-th entry is given by $``\nabla_i d\cdot\nabla_j d"$ for $1\leq i,j\leq 3$.

There is a vast literature on the hydrodynamics of the liquid crystal system. For background we list  a few names, with no intention to be complete:  \cite{Cal, CC, CC1, CDLL, CR, CRW,  DQSr, EK, JT,  HKL, Kin, LL, LL2, LL1, Liu, SL, WXL}.  Particularly, in \cite{WXL}, the general Ericksen-Leslie system with constant density is studied.  In the work, the global regularity and long-time behavior of solutions are obtained with the assumption that the viscosity coefficient is sufficiently large (3D). With the Parodi's relation, the authors established the global  well-posedness and Lyapunov stability near local energy minimizers. The authors also discussed the connection between
Parodi's relation and the linear stability. On the other hand, considering the density is not constant, a relatively full model for the dynamic of Smectic-A liquid crystals is studied in \cite{Lc2}. In this model, several terms in $\sigma$ are assumed to be zero, and the term $-\Omega d+\frac{\lambda_2}{\lambda_1}Ad$ does not appear in equation (\ref{eqL3}). The author proved the existence of global classical solutions in both two and three dimensional cases. In 2D, no additional assumption is needed; while in 3D, the flow viscosity coefficient is assumed to be sufficiently large. Moreover, a ``weak-strong" type of uniqueness result, a long-time behavior and stability of solutions are obtained. 
Also, with a non-constant density, a simplified Ericksen-Leslie model is studied in \cite{DQSr}. The authors indicated that a regular solution exists globally in 2D with general data, while exists globally in 3D with small initial data or for short time with general data.

 In the present paper the consideration is given to the full Ericksen-Leslie model (\ref{eqL1})-(\ref{eqL4}) with non-constant density and under the Parodi's relation (\ref{eq:Par}).  We establish that: in 2D, there exists a global regular solution to the system with general data (no extra condition); in 3D, there exists a global regular solution with small initial data or a local (short time) regular solution with general data.  We also show that a ``weak-strong" type of uniqueness result holds with certain assumption on the weak solution. Namely, if there exist a regular solutions satisfying a
 higher order energy estimate and a weak solution satisfying the basic energy estimate and two auxiliary estimates, they must be identical.

For the simplified liquid crystal model with constant density (cf. \cite{LL}) or non-constant density (cf. \cite{DQSr}), and for the general Ericksen-Leslie system (\ref{eqL1})-(\ref{eqL4}) with the artificial assumption $\lambda_2=0$ (cf. \cite{LL1}), the transport equation of $d$ satisfies a certain type of maximum principle. In the present paper, for the general model (\ref{eqL1})-(\ref{eqL4}) with $\lambda_2\neq 0$, the stretching effect causes the loss of maximum principle  for $d$ (cf. \cite{WXL}).  However, in the analysis of sequel, the estimate
$$
d\in L^\infty(0,T;L^\infty)
$$
turns out to be essential in the derivation of higher-order energy estimates and hence to assure that the stress term $\nabla\cdot\sigma$ can be handled successfully. Therefore, as in \cite{WXL}, we consider the periodic boundary conditions which help to avoid the difficulties from boundary terms when deriving the higher order energy estimates. As such, we restrict ourselves to the following boundary conditions:
\bg\label{bd}
u(x+e_i,t)=u(x,t),  \ \ d(x+e_i, t)=d(x,t), \ \ \ \mbox { for } (x,t)\in\partial Q\times \R^+
\ed
with $Q$ being a unit square in $\R^2$ or cube in $\R^3$. 
The full Ericksen-Leslie model also satisfies the initial conditions:
\bg\label{initden}
\rho(x,0) =\rho_0(x), \ \ \ 0<M_1\leq\rho_0(x)\leq M_2,
\ed
\bg\label{initu}
u(x,0) =u_0(x), \qquad \nabla\cdot u_0 =0, \qquad \mbox { and }  \qquad  d(x,0) =d_0(x).
\ed

 In the rest of the introduction we describe  our main results: \\

%\noindent \emph{Regularity in 2 dimensions:}
\begin{Theorem}\label{Mthm1}(2D) Suppose that $Q\subset\R^2$ is a unit square. Let $\rho_0$, $u_0$ and $d_0$ satisfy (\ref{initden}) and (\ref{initu}).
Suppose that $\rho_0\in C^1$, $u_0\in H_p^2$ and $d_0\in H_p^{3}$. Then, system (\ref{eqL1})-(\ref{eqL4}) has a global classical solution $(\rho,u,d)$, that is, for all $T>0$ and some $\alpha\in(0,1)$
\begin{equation}\label{reg}\begin{split}
u\in C^{1+\alpha/2,2+\alpha}((0,T)\times Q),\\
\nabla p\in C^{\alpha/2,\alpha}((0,T)\times Q),\\
d\in C^{1+\alpha/2,2+\alpha}((0,T)\times Q),\\
\rho\in C^1((0,T)\times Q).
\end{split}\end{equation}
\end{Theorem}

The notations $H_p^2$ and $H_p^3$ will be introduced in Section \ref{sec:pre}.\\
%The regularity for the flow of the nematic liquid crystals in dimension 2 turns out to be not too difficult. We first establish the Ladyzhenskaya energy estimate \eqref{phiest} and \eqref{L^2H^2}, similar to that in \cite{LL}. Then we apply the regularity result for transport equations in \cite{AKM} to obtain the H\"{o}lder continuity of the fluid density. Therefore
%Theorem \ref{Mthm1} follows from the $L^p$ estimates and H\"older estimates in \cite{LS} and a more or less standard bootstrapping between the three equations in the
%system \eqref{LCD}.\\

%similar to the argument to prove the regularity of the solutions to the pure fluid systems in \cite{AKM}.

Provided we have sufficiently small data or we work with sufficiently short time, we also obtain the regularity in three dimensional case. %Given more restriction on data, the regularity can be obtained up to the boundary.\\

%\noindent \emph{Regularity  in  3 dimension  with small data or short time:}
\begin{Theorem}\label{Mthm2} (3D) Suppose that $Q\subset\R^3$ is a unit cube. Let $\rho_0$, $u_0$ and $d_0$ satisfy (\ref{initden}) and (\ref{initu}). Assume that $\rho_0\in C^1$, $u_0\in H_p^2$ and $d_0\in H_p^{3}$.   Then\\
$\mathbf{1}$. There is a positive small number $\e_0$ such that if
\begin{equation}\label{smalli}
\rho_0\|u_0\|_{H_p^1}^2+\|d_0\|_{H_p^1}^2+\|\Delta d_0-f(d_0)\|_{L_p^2}^2\leq\epsilon_0,
\end{equation}
then system (\ref{eqL1})-(\ref{eqL4}) has a classical solution $(\rho,u,d)$ in the time period $(0, T)$, for all $T>0$. That is, (\ref{reg}) holds for some $\alpha\in(0,1)$.

\noindent $\mathbf{2}$. For general data, there exists a positive number $\delta_0 = \delta_0(\rho_0, u_0, d_0)$ such that (\ref{reg}) holds in the interval $(0,T)$ for some small $T\leq\delta_0$.
\end{Theorem}

%We obtain interior regularity with a relatively weak conditions on the initial data. For more regular data, we are able to obtain solutions which are regular up to the boundary. Indeed, more regularity on the initial data $d_0$ implies more regularity on the boundary due to the second condition in (\ref{bd}). Thus, yielding the regularity up to boundary, see \cite{Kr1}.

%\begin{Corollary}\label{Cor1}
%Suppose $(\rho_0,u_0,d_0)$ satisfies the conditions in Theorem \ref{Mthm1}. In addition, $\rho_0\in C^1(\bar Q)$,  $u_0\in C^{2+\alpha}(\bar Q)$ and $d_0\in C^{2+\alpha}(\bar Q)$. Then the solution $(\rho,u,d)$ is regular up to boundary, that is, the conclusions in Theorem \ref{Mthm1} hold in$[0,T]\times\bar Q$.
%\end{Corollary}

\begin{Remark}
It is pointed out in \cite{Liu, WXLn} that the large viscosity assumption is not equivalent to the small initial data assumption for the Ericksen-Leslie system (\ref{eqL1})-(\ref{eqL4}), due to its much more complicated structure. Thus it is particularly interesting to investigate the regularity of solutions to the Ericksen-Leslie system under the small initial data assumption.
\end{Remark}

\begin{Remark}
In contrast to the initial condition in \cite{DQSr} (Theorem 1.3) where $(u_0,\nabla d_0)\in H^1\times H^1$, we require a higher order condition on the initial data here, that is $(u_0, \nabla d_0)\in H_p^2\times H_p^2$. The reason is that, in the full system, the term $\nabla \cdot \sigma$ in (\ref{eqL2}) contains the high order term $\nabla\cdot(Nd)$ which is not in the simplified system in \cite{DQSr}. To deal with the high order term $\nabla\cdot(Nd)$, we need a higher order energy estimate compared to the case of the simplified system. But to obtain the higher order energy estimate, the smallness assumption on  $(u_0,\nabla d_0)$ in $H_p^1\times H_p^1$ is sufficient. %{\color{red}double check it really needs one order higher estimate.}
\end{Remark}

%The proof of the regularity of solutions to system \eqref{LCD} in three dimension takes the same approach as in two dimension but is much more complicated. First in contrast to the case of two dimension, we only get the Ladyzhenskaya energy estimates when either the initial data is small in the sense as described in (\ref{smalli}) or $T$ is small. Our calculations and estimates are based on those in \cite{LL}, with interesting modifications. We use ideas of \cite{LL} making it work in a rather different way. We keep the potentially small terms $\|u\|_{L^2}$ and $\|\nabla d\|_{L^2}$ instead of throwing them away. This gives a more unified way to derive the Ladyzhenskaya energy estimates in the cases of small data or for short time.

%After having the Ladyzhenskaya energy estimates, in contrast to the two dimensional case, we do not have the H\"{o}lder continuity for the fluid density. Instead we observe that we have small oscillations of the density over small balls in $\Omega\times[0,T]$ provided that either the initial data is small or for short time, see \cite{DQSr} for a detailed proof. This turns out to be enough to carry out the frozen coefficient method to improve the regularity of the fluid velocity. We refer the reader to \cite{LS1} for a reference of the frozen coefficient method. \\

The basic idea of the regularity proof is to get some high order energy estimates (Ladyzhenskaya method) as described in \cite{LL, DQSr}. To fully utilize the smallness assumption (\ref{smalli}), we follow the principle in \cite{DQSr} by keeping the potentially small terms $\|u\|_{L^2}$ and $\|\nabla d\|_{L^2}$ instead of throwing them away. However, we point out that the situation for the full system in the present work is much more complicated than that of the simplified system in \cite{DQSr}, because the velocity field equation (\ref{eqL2}) and director field equation (\ref{eqL3}) are coupled not only through $\nabla \cdot (\nabla d\otimes\nabla d)$ but also through $\nabla\cdot\sigma$ and $-\Omega d+\frac{\lambda_2}{\lambda_1} Ad$, while $\sigma$ consists of six non-zero terms.

%\begin{Remark}
%The interior regularity in Theorem \ref{Mthm1} and Theorem \ref{Mthm2} is obtained by bootstrapping argument.
%\end{Remark}

In the following, we state the uniqueness of solutions,
\begin{Theorem}\label{unique}
Let $(\rho,u,d)$ be a solution to system (\ref{eqL1})-(\ref{eqL4}) and (\ref{bd})-(\ref{initu}) obtained in Theorem \ref{Mthm1} for two dimensional case or in Theorem \ref{Mthm2} for three dimensional case.
%We suppose the additional conditions on the data:  $\rho_0\in C^1(\bar\Omega)$, $u_0\in C^{2+\alpha}(\Omega)$ and $d_0\in C^{2+\alpha}(\Omega)$.
Let $(\bar\rho,\bar u,\bar d)$ be a weak solution to system (\ref{eqL1})-(\ref{eqL4}) with (\ref{bd})-(\ref{initu}) satisfying the following energy inequalities:
\begin{align}\label{weakeng}
&\int_{Q}|\bar\rho|^2 dx \leq \int_{Q}|\rho_0|^2 dx,\\
&\int_{Q}\frac 1 2\bar\rho|\bar u|^2+\frac 1 2|\nabla\bar d|^2+
\bar\rho F(\bar d)dx+\int_0^T\int_{Q}\mu_4|\nabla\bar u|^2-\frac{1}{\lambda_1}|\Delta\bar d-f(\bar d)|^2dxdt\\
&\leq \int_{Q}\frac 1 2\rho_0|u_0|^2+\frac 1 2|\nabla d_0|^2+\rho_0F(d_0)dx.\notag
\end{align}
In addition, $|\bar d|$ is bounded.
Then $(\rho,u,d)\equiv (\bar\rho,\bar u,\bar d)$.
\end{Theorem}

The uniqueness result is achieved by a standard approach where one establishes a Gronwall's inequality for the difference of the two solutions. A uniqueness proof for the simplified system was presented in \cite{DQSr}. For the full system the process will be similar. Since the computation work is huge and tedious, we give the proof of Theorem \ref{unique} in Appendix (Section \ref{sec:uni}).

The rest of the paper is organized as: in Section \ref{sec:pre} we introduce some notations that shall be used throughout the paper and the basic energy law governing the full system (\ref{eqL1})-(\ref{eqL4}); in Section \ref{sec:proof} we prove Theorems \ref{Mthm1} and \ref{Mthm2} by several steps; in Appendix we devote to proving Theorem \ref{unique}.

%8888888888888888888888888888888888888888888888888888888888888888888888888

\bigskip

\section{Preliminary} 
\label{sec:pre}

\subsection{Notations}
We adopt the standard functional settings and notations for periodic  problems (cf. \cite{Tem1}) in the following:
\begin{align}\notag
&H^m_p(Q)=\left\{u\in H^m(\R^n,\R^n)|u(x+e_i)=u(x)\right\}\\
&\dot H^m_p=H^m_p(Q)\cap \left\{\int_Qu(x)dx=0\right\}\notag\\
&H=\left\{u\in L^2_p(Q)|\nabla\cdot u=0\right\}  \ \ \ \mbox { with } L^2_p(Q)=H^0_p(Q)\notag\\
&V=\left\{u\in \dot H^1_p(Q)|\nabla\cdot u=0\right\}\notag\\
&V'=\mbox { the dual space of } V.\notag
\end{align}
At certain places in the paper, we also denote the space of scalar functions by $H^m_p(Q)$.

The inner product on $L^2_p(Q)$ and $H$ is denoted by $(\cdot,\cdot)$ and the associated norm by $\|\cdot\|$. For simplicity, the space $H^m_p(Q)$ is denoted by $H^m_p$. The inner product on $H^m_p$ is defined as $<u, v>_m=\Sigma_{|k|=0}^m(D^ku,D^kv)$ with $k=(k_1,\ldots,k_n)$ being the multi-index of length $|k|=\Sigma_{i=1}^nk_i$ and $D^k=\partial_{x_1}^{k_1}\ldots\partial_{x_1}^{k_1}$.

\subsection{Definition of weak solution}
The weak formulation of the problem is given as follows.
\begin{Definition}\label{def:weak}
The triplet $(\rho, u,d)$ is called a weak solution to the system (\ref{eqL1})-(\ref{eqL4}) in $Q_T=Q\times(0,T)$ subject to the boundary and initial conditions (\ref{bd})-(\ref{initu}), if it satisfies 
\begin{align}\notag
&0<M_1\leq\rho\leq M_2 \\
&u\in L^\infty(0,T;H)\cap L^2(0,T;V)\notag\\
&d\in L^\infty(0,T; H^1_p\cap L^\infty_p)\cap L^2(0,T;H^2_p)\notag
\end{align}
and moreover, for any smooth function $\psi(t)$ with $\psi(T)=0$ and $\phi(x)\in H^1_p$, the following integral equations hold:
\begin{equation}\notag%\label{int:u}
%\begin{split}
\int_0^T(\rho,\psi_t\phi)dt-\int_0^T(\rho u,\psi\nabla\phi)dt
=(\rho_0,\phi)\psi(0),
%\end{split}
\end{equation}
\begin{equation}\notag%\label{int:u}
\begin{split}
&\int_0^T(u,\psi_t\phi)dt-\int_0^T(u\cdot\nabla u,\psi\phi)dt\\
=&(u_0,\phi)\psi(0)-\int_0^T(\nabla d\otimes\nabla d,\psi\nabla\phi)dt+\int_0^T(\sigma,\psi\nabla\phi)dt,
\end{split}
\end{equation}
\begin{equation}\notag%\label{int:u}
\begin{split}
&\int_0^T(d,\psi_t\phi)dt-\int_0^T(u\cdot\nabla d,\psi\phi)dt+\int_0^T(\Omega d,\psi\phi)dt-\frac{\lambda_2}{\lambda_1}\int_0^T(Ad,\psi\phi)dt\\
=&(d_0,\phi)\psi(0)+\frac{1}{\lambda_1}\int_0^T(\Delta d-f(d),\psi\phi)dt.
\end{split}
\end{equation}
\end{Definition}

\subsection{Basic energy law}

The total energy of the full Ericksen-Leslie system (\ref{eqL1})-(\ref{eqL4}) is given by 
\bg\label{id:en}
\mathcal E(t)=\int_Q\frac{1}{2}\rho |u|^2+\frac{1}{2}|\nabla d|^2+\rho \mathcal F(d)dx
\ed
which consists of kinetic and potential energies. Formally, a smooth solution $(\rho,u,d)$ satisfies (cf. \cite{LL1, Lc2})
\begin{align}\label{eq:en}
\frac{d}{dt}\mathcal E(t)=&-\int_Q\mu_1|d^TAd|^2+\mu_4|\nabla u|^2+(\mu_5+\mu_6)|Ad|^2dx\\
&+\lambda_1\|N\|^2+(\lambda_2-\mu_2-\mu_3)(N,Ad).\notag
\end{align}
To guarantee the dissipation of the director field, it is assumed that (cf. \cite{Er2,Le1})
\begin{equation}\label{ass:para}
\begin{cases}
&\lambda_1<0\\
&\mu_1\geq 0, \ \ \ \mu_4>0,\\
&\mu_5+\mu_6\geq 0.
\end{cases}
\end{equation}

With the Parodi's relation (\ref{eq:Par}) we have the following basic energy law.

\begin{Lemma}\label{le:basic-en}
Suppose that (\ref{re:con}), (\ref{eq:Par}) and (\ref{ass:para}) are satisfied. In addition, we assume 
\bg\label{ass:para2}
\frac{\lambda_2^2}{-\lambda_1}\leq\mu_5+\mu_6.
\ed
Then the total energy $\mathcal E(t)$ for smooth solution satisfies 
\begin{align}\label{en:dec}
\frac{d}{dt}\mathcal E(t)&=-\int_Q\mu_1|d^TAd|^2+\mu_4|\nabla u|^2dx+\frac{1}{\lambda_1}\|\Delta d-f(d)\|^2\\
&-\left(\mu_5+\mu_6+\frac{\lambda_2^2}{\lambda_1}\right)\|Ad\|^2\leq 0.\notag
\end{align}
\end{Lemma}
The proof of Lemma \ref{le:basic-en} is identical to the proof of Lemma 2.1 in \cite{WXL} and thus it is omitted here.

\bigskip

\section{Regular solutions to the full Ericksen-Lesie System}
\label{sec:proof}

\subsection{Galerkin Approximate Solutions}
We construct a sequence of Galerkin approximating solutions that satisfy both of the basic energy estimate and a higher order energy estimate. The higher order energy estimate obtained through Ladyzhenskaya method yields a subsequence that will converge to the classical solution.\\
We first introduce the functional settings of Galerkin approximation method.

Let
$$
\Hh(Q)= \mbox { closure of } \{f\in C^\infty_0(Q,\mathbb{R}^n): \nabla\cdot f = 0\}   \mbox { in } H.
$$
Let
$\left\{\phi_i\right\}_{i=1}^\infty$ be the unit eigenvectors of the Stokes problem in the periodic case with zero mean:
\bg\notag
-\triangle\phi_i+\nabla \pi_i=k_i\phi_i \mbox { in } Q, \ \ \int_Q\phi_i(x)dx=0
\ed
with $\pi_i\in L^2$ and $0<k_1\leq k_2\leq\ldots$ being eigenvalues.
It is known that $\phi_i$ are smooth and $\left\{\phi_i\right\}_{i=1}^\infty$ forms an orthogonal basis of $\Hh$ (see \cite{Tem1}). Let
$$
P_m: \Hh\to \Hh_m=span\left\{ \phi_1,...\phi_m\right\}
$$
be the orthonormal projection. We seek approximate solutions $(\rho^m, u^m, d^m)$ with $u^m\in \Hh_m$, satisfy the following equations:
\bg\label{aden}
\rho_t^m+u^m\cdot\nabla\rho^m=0,
\ed
\bg\label{aNSE}
P_m(\rho^m\frac{\partial}{\partial t}u^m)=P_m(\triangle u^m-\rho^mu^m\cdot\nabla u^m-\nabla\cdot(\nabla d^m\otimes\nabla d^m)-\nabla\cdot \sigma^m)
\ed
\bg\label{ad}
d_t^m+u^m\cdot\nabla d^m-\Omega^md^m+\frac{\lambda_2}{\lambda_1}A^md^m=-\frac{1}{\lambda_1}(\Delta d^m-f(d^m))
\ed
with
\begin{equation}\notag
\begin{split}
\Omega^m=&\frac{1}{2}(\nabla u^m-\nabla^Tu^m),\ \ \ A^m=\frac{1}{2}(\nabla u^m+\nabla^Tu^m),\\
N^m=&\partial_td^m+(u^m\cdot\nabla)d^m+\Omega^md^m,\\
\sigma^m=&\mu_1((d^m)^TA^md^m)d^m\otimes d^m+\mu_2N^m\otimes d^m+\mu_3d^m\otimes N^m+\mu_4A^m\\
&+\mu_5A^md^m\otimes d^m+\mu_6 d^m\otimes A^md^m,
\end{split}
\end{equation}
and with the initial and boundary conditions
\bg\notag%\label{auinit}
\rho^m(x,0)=\rho_0(x), \quad u^m(x,0)=P_m u_0(x), \quad
d^m(x,0)=d_0(x),
\ed
\bg\notag%\label{aubd}
u^m(x+e_i,t)=u^m(x,t), \quad d^m(x+e_i,t)=d^m(x,t).
\ed

Let
\bg\notag%\label{au}
u^m(x,t)=\sum_{l=1}^mg_i^m(t)\phi_i(x),
\ed
with  $g_i^m(t)\in C^1[0,T]$. Hence (\ref{aNSE}) is equivalent to the following system of ordinary differential equations:
\bg\label{ODE}
\sum_{i=1}^m\Aa_i^{mj}(t)\frac{d}{dt}g_i^m(t)=-\sum_{i,k}^m\Bb_{ik}^{mj}(t)g_i^m(t)g_k^m(t)
-\sum_{i=1}^m\Cc_i^jg_i^m(t)+\Dd^{mj}(t),
\ed
for $j=1,2,..., m$, where
\begin{equation}\notag%\label{cff}
\begin{cases}
\Aa_i^{mj}(t)=\int_\Omega\rho^m(t)\phi_i(x)\phi_j(x)dx,\\
\Bb_{ik}^{mj}(t)=\int_\Omega\rho^m(t)(\phi_i(x)\cdot\nabla\phi_k(x))\phi_j(x)dx,\\
\Cc_i^{j}=\int_\Omega\nabla\phi_i(x)\cdot\nabla\phi_j(x)dx,\\
\Dd^{mj}(t)=\int_\Omega\sum_{k,l}(\frac{\partial}{\partial x_k}d^m\cdot\frac{\partial}{\partial x_l}d^m+\sigma^m_{kl})\frac{\partial}{\partial x_l}\phi_j^k(x)dx.
\end{cases}\end{equation}
Here $\phi_j^k(x)$ is the $k$-th component of the vector $\phi_j(x)$.  And
$$
u^m(\cdot, 0) = \sum_{i=1}^mg_i^m(0)\phi_i(x),\quad \text{where} \quad
g_i^m(0)=\int_\Omega u_0(x)\phi_i(x)dx.
$$

\begin{Lemma}\label{le:weak}
There exists a solution $(\rho^m,u^m,d^m)$ to the problem (\ref{aden})-(\ref{ad}) in $Q_T=Q\times[0,T]$, for any $T\in(0,\infty)$, satisfying
$$
M_1 \leq \rho^m \leq M_2
$$
$$
u^m\in L^\infty(0,T;H)\cap L^2(0,T;V)
$$
$$
d^m\in L^\infty(0,T; H^1_p\cap L^\infty_p)\cap L^2(0,T;H^2_p).
$$
Moreover, $(\rho^m,u^m,d^m)$ is smooth in the interior of $Q_T$ and satisfies the basic energy equality,
\begin{equation}\label{energy}
\begin{split}
\frac{d}{dt}\mathcal E^m(t)=&-\int_Q\mu_1|(d^m)^TA^md^m|^2+\mu_4|\nabla u^m|^2dx+\frac{1}{\lambda_1}\|\Delta d^m-f(d^m)\|^2\\
&-\left(\mu_5+\mu_6+\frac{\lambda_2^2}{\lambda_1}\right)\|A^md^m\|^2\leq 0,
\end{split}
\end{equation}
with 
\bg\notag
\mathcal E^m(t)=\int_Q\frac{1}{2}\rho^m |u^m|^2+\frac{1}{2}|\nabla d^m|^2+\rho^m \mathcal F(d^m)dx.
\ed
\end{Lemma}
\pf
The proof of existence of weak solutions is based on an application of the Leray-Schauder fixed point theorem. Let $v^m=\sum_{i=1}^m h^m_i \phi_i \in C^1(0, T; \Hh_m)$. For each $m$ let $\rho^m$ be a solution to
\begin{equation}\label{transport}
\rho + v^m\cdot \nabla\rho = 0
\end{equation}
with initial condition $\rho(\cdot, 0) = \rho_0$. Let $d^m$ be a solution to
$$
d_t+v^m\cdot\nabla d-\tilde\Omega^md+\frac{\lambda_2}{\lambda_1}\tilde A^md=-\frac{1}{\lambda_1}(\Delta d-f(d))
$$
$$
\tilde\Omega^m=\frac{1}{2}(\nabla v^m-\nabla^Tv^m),\ \ \ \tilde A^m=\frac{1}{2}(\nabla v^m+\nabla^Tv^m)
$$
with initial condition $d(\cdot, 0) = d_0$ and boundary condition $d(x+e_i,t)=d(x,t)$. The reason the transport equation \eqref{transport} is solvable for $v^m\in C^1(0, T; \Hh_m)$ is due to the regularity of the eigenfunctions of the Stokes operators (cf. \cite{Tem} \cite{La}). Let
$u^m = \sum_{i=1}^m g^m_i \phi_i\in C^1(0, T; \Hh_m)$ be the solution of the system of linear equations
$$
\sum_{i=1}^m\Aa_i^{mj}(t)\frac{d}{dt}g_i^m(t)=-\sum_{i,k}^m\Bb_{ik}^{mj}(t)h_i^m(t)g_k^m(t)
-\sum_{i=1}^m\Cc_i^jg_i^m(t)+\Dd^{mj}(t).
$$
This system of linear equations is solvable because the eigenvalues of matrix of the coefficients $\Aa_i^{mj}(t)$ are bounded from below, since
\begin{equation}\label{inverse}
\Aa_i^{mj}\gamma_i\gamma_j = \int_\Omega \rho |\psi|^2 dx \geq M_1\int_\Omega |\psi|^2 dx\quad \text{where $\psi = \sum_{i=1}^m \gamma_i\phi_i$}.
\end{equation}
Thus we constructed a mapping $\Mm$ with $\Mm(v^m)=u^m$.
The energy estimate (\ref{energy}) will be obtained similarly as in Lemma \ref{le:basic-en} and it plays the key to allow one to apply Leray-Schauder fixed point theorem for $\Mm$.\\

\cbdu

\bigskip

\subsection{Higher order energy estimates}
\label{sec:higher-energy}
As a first step, we use the Ladyzhenskaya energy method \cite{La} to show that $u^m \in L^\infty(0, T; V)\bigcap L^2(0, T; H^2_p)$ and $d^m\in L^\infty(0, T; H_p^2)$ $\bigcap$ $L^2(0, T; H_p^3)$, provided $u_0\in H_p^1$ and $d_0\in H_p^2$. In 2D, these estimates will be obtained with general initial data; while in 3D, they will be obtained under the assumption of small initial data.  We then pass to the limit for the Galerkin approximating solutions $(\rho^m,u^m,d^m)$ to yield weak solutions for the system (\ref{eqL1})-(\ref{eqL4}). The key inequalities used often in this paper are the following Gagliardo-Nirenberg inequality (cf. \cite{Ev}) :

\begin{Lemma}\label{le:GN} (Gagliardo-Nirenberg) If $\Omega$ is a smooth bounded domain in $\mathbb{R}^n$, then
\bg\notag%\label{GNd2+}
\|u\|_{L^4(\Omega)}^4\leq C\|u\|_{L^2(\Omega)}^2(\|\nabla u\|_{L^2(\Omega)}^2 + \|u\|_{L^2(\Omega)}^2), \qquad n=2,
\ed
\bg\notag%\label{GNd3+}
\|u\|_{L^4(\Omega)}^4\leq C\|u\|_{L^2(\Omega)}(\|\nabla u\|_{L^2(\Omega)}^2 + \|u\|_{L^2(\Omega)}^2)^\frac 32, \qquad n=3.
\ed
%Moreover, if $u|_{\partial\Omega} = 0$, then
%\bg\label{GNd2}
%\|u\|_{L^4(\Omega)}^4\leq C\|u\|_{L^2(\Omega)}^2\|\nabla u\|_{L^2(\Omega)}^2, \qquad n=2,
%\ed
%\bg\label{GNd3}
%\|u\|_{L^4(\Omega)}^4\leq C\|u\|_{L^2(\Omega)}\|\nabla u\|_{L^2(\Omega)}^3, \qquad n=3.
%\ed
\end{Lemma}

%888888888888888888888
%\subsection{The $L^\infty(H^1)$ and $L^2(H^2)$ estimates of the velocity}

By a similar strategy as in \cite{LL, DQSr}, we will establish uniform estimates for the following energy quantity
\bg\notag%\label{phi}
\Phi_m^2(t)=\|\sqrt{\rho^m}\nabla u^m\|_{L^2}^2+\|\triangle d^m-f(d^m)\|_{L^2}^2.
\ed
%Then to show that $u^m \in L^\infty( H^1)\bigcap L^2(H^2)$ and $d^m\in L^\infty( H^2)\bigcap L^2(H^3)$ we calculate $\frac d{dt} \Phi_m^2$.  Namely,

%In 3D we have: 

\begin{Lemma}\label{le:phi} (3D small data assumption.) Let $Q\subset\mathbb{R}^3$ and   the initial data $(\rho_0,u_0,d_0)$  satisfy (\ref{initden})-(\ref{initu}). Suppose  that %$u_0\in H_p^1$, $d_0\in H_p^2$ and
 $\Phi^2(0)=\|\sqrt{\rho_0}\nabla u_0\|_{2}^2+\|\Delta d_0-f(d_0)\|_{2}^2<\infty$. There is $\epsilon_0>0$ such that if
\begin{equation}\label{ass:small}
(\frac{M_2}{\mu_4}+1)\rho_0\|u_0\|_{H^1}^2 + \|d_0\|_{H^1}^2+\|\Delta d_0-f(d_0)\|_2^2 \leq \epsilon_0,
\end{equation}
then the approximating solutions $(\rho^m,u^m,d^m)$ obtained in Lemma \ref{le:weak}  satisfy
\begin{equation}\label{phiest}
\|\sqrt{\rho^m}\nabla u^m\|_{L^2}^2 + \|\Delta d^m-f(d^m)\|_{L^2}^2 \leq C(\Phi^2(0)+C)
\end{equation}
for all $t\in [0, T]$ and
\bg\label{L2H2}
\int_0^T\|\Delta u^m\|_{2}^2 +\|\nabla(\Delta d^m-f(d^m))\|_{2}^2+\mu_1\|d^{mT}\nabla A^md^m\|_2^2dt \leq C
\ed
where the constants $C$ depend only on initial data, $Q$, $M_1$, $M_2$ and the physical coefficients in system (\ref{eqL1})-(\ref{eqL4}).
\end{Lemma}
\pf
Note that the approximating solutions $(\rho^m, u^m, d^m)$ obtained in Lemma \ref{le:weak} also satisfy equations (\ref{eqL1})-(\ref{eqL3}) point wisely.
To simplify the notation, through the proof we drop the approximating index $m$ and denote $(\rho, u, d)$ as the Galerkin approximating solutions $(\rho^m, u^m, d^m)$. We also denote $\Phi_m(t)$ by $\Phi(t)$.
%\begin{align}\notag
%\frac{1}{2}\frac{d}{dt}\Phi_m^2(t)&=\int_{\Omega}\nabla u^m\cdot\nabla u_t^mdx+\int_{\Omega}(\Delta d^m-f(d^m))\cdot(\Delta d_t^m-f'(d^m)d_t^m)dx\\
%&= - \int_{\Omega}\rho^m|u^m_t|^2dx - \int_\Omega |\nabla\triangle d^m|^2dx\notag\\
%&+\int_{\Omega}\nabla\triangle d^m\cdot\nabla(u^m\cdot\nabla d^m)-\triangle d^m\cdot\triangle(f(d^m))dx\notag\\
%&+\int_{\Omega}-\rho^m(u^m\cdot\nabla u^m)u_t^m-u_t^m\nabla d^m\triangle d^mdx\notag\\
%&= - \int_{\Omega}\rho^m|u^m_t|^2dx - \int_\Omega |\nabla\triangle d^m|^2dx\notag\\
%& + I + II + II + IV,\notag
%\end{align}
By equation (\ref{eqL1}) and (\ref{eqL4}), applying the periodic boundary conditions (\ref{bd}) and integration by parts yields 
\begin{equation}\label{eq:phi0}
\begin{split}
&\frac{1}{2}\frac{d}{dt}\Phi^2(t)\\
=&\int_{Q}\frac{1}{2}\rho_t|\nabla u|^2+\rho\nabla u\nabla u_tdx+\int_{Q}(\Delta d-f(d))\cdot(\Delta d_t-f'(d)d_t)dx\\
=&\int_{Q}\frac{1}{2}\rho (u\cdot\nabla)|\nabla u|^2-\rho u_t\Delta udx+\int_{Q}(\Delta d-f(d))\Delta d_tdx\\
&-\int_{Q}(\Delta d-f(d))f'(d)d_tdx.
\end{split}
\end{equation}

It follows from equation (\ref{eqL2})
\begin{equation}\label{eq:en-A}
\begin{split}
&\int_{Q}\frac{1}{2}\rho (u\cdot\nabla)|\nabla u|^2-\rho u_t\Delta udx\\
=&\int_{Q}\frac{1}{2}\rho (u\cdot\nabla)|\nabla u|^2+\rho (u\cdot\nabla)u\Delta u+\Delta u\nabla d\Delta d-\Delta u\nabla\cdot \sigma dx\\
=&-\mu_4\int_Q|\Delta u|^2dx+\int_{Q}2\rho (u\cdot\nabla)u\Delta u+\Delta u\nabla d\Delta d-\Delta u\nabla\cdot \bar\sigma dx\\
=&-\mu_4\int_Q|\Delta u|^2dx+I_1+I_2+I_3
\end{split}
\end{equation}
with $\bar\sigma=\sigma-\mu_4A$. 
Acting $Laplacian$ on equation (\ref{eqL3}) gives
\begin{equation}\notag
\Delta d_t+\Delta(u\cdot\nabla d)-\Delta(\Omega d)+\frac{\lambda_2}{\lambda_1}\Delta(Ad)
=-\frac{1}{\lambda_1}\Delta(\Delta d-f(d)).
\end{equation}
Thus,
\begin{equation}\label{eq:en-dt}
\begin{split}
&\int_{\Omega}(\Delta d-f(d))\Delta d_tdx\\
=&\frac{1}{\lambda_1}\int_\Omega|\nabla(\Delta d-f(d))|^2dx-\int_\Omega(\Delta d-f(d))\Delta(u\cdot\nabla d)dx\\
&+\int_\Omega(\Delta d-f(d))\Delta(\Omega d)dx-\frac{\lambda_2}{\lambda_1}\int_\Omega(\Delta d-f(d))\Delta(Ad)x\\
=&\frac{1}{\lambda_1}\int_\Omega|\nabla(\Delta d-f(d))|^2dx+I_4+I_5+I_6.
\end{split}
\end{equation}
We also have the following by equation (\ref{eqL3}) 
\begin{equation}\label{eq:en-f}
\begin{split}
&-\int_{\Omega}(\Delta d-f(d))f'(d)d_tdx=\frac{1}{\lambda_1}\int_{\Omega}f'(d)|\Delta d-f(d)|^2dx\\
&+(\Delta d-f(d), f'(d)(u\cdot\nabla)d)-(\Delta d-f(d), f'(d)(\Omega d-\frac{\lambda_2}{\lambda_1}Ad))\\
&=I_7+I_8+I_9.
\end{split}
\end{equation}
Combining (\ref{eq:phi0}), (\ref{eq:en-A}), (\ref{eq:en-dt}) and (\ref{eq:en-f}) yields
\begin{align}\label{eq:phi1}
\frac{1}{2}\frac{d}{dt}\Phi^2(t)=&-\mu_4\int_{Q}|\Delta u|^2dx+\frac{1}{\lambda_1}\int_Q|\nabla(\Delta d-f(d))|^2dx\\
&+I_1+I_2+I_3+I_4+I_5+I_6+I_7+I_8+I_9\notag,
\end{align}
with $I_k$ being defined as in the above equations for $k=1, \ldots, 9$. In the following, we shall estimate these $I_k$ term by term. 

By H\"older inequality and Gagliardo-Nirenberg inequality in Lemma \ref{le:GN}, it follows that 
\begin{equation}\label{esti:i1}
\begin{split}
|I_1|&\leq M_2\int_\Omega|u||\nabla u||\Delta u|dx\\
&\leq \epsilon\|\Delta u\|_2^2+C\|u\|_4^2\|\nabla u\|_4^2\\
&\leq \epsilon\|\Delta u\|_2^2+C\|u\|_2^{1/2}\|\nabla u\|_2^2\|\Delta u\|_2^{3/2}\\
&\leq \epsilon\|\Delta u\|_2^2+C\|u\|_2^2\|\nabla u\|_2^8\\
&\leq \epsilon\|\Delta u\|_2^2+CM_1^{-4}\|u\|_2^2\|\sqrt{\rho}\nabla u\|_2^8.
\end{split}
\end{equation}
By Agmon's inequality in 3D, we have
\begin{equation}\label{Ag}
\|d\|_\infty\leq C\|\nabla d\|_2^{1/2}\|\Delta d\|_2^{1/2}. 
\end{equation}
Thus we infer that
\begin{equation}\label{est:d-h2}
\|\Delta d\|_2\leq\|\Delta d-f(d)\|_2+\|f(d)\|_2\leq\|\Delta d-f(d)\|_2+C,
\end{equation}
\begin{equation}\label{esti:d-h3}
\begin{split}
\|\nabla\Delta d\|_2&\leq \|\nabla(\Delta d-f(d))\|_2+\|\nabla f(d)\|_2\\
&\leq \|\nabla(\Delta d-f(d))\|_2+\|f'(d)\|_\infty\|\nabla d\|_2\notag\\
&\leq \|\nabla(\Delta d-f(d))\|_2+(1+\|d\|_\infty^2)\|\nabla d\|_2\notag\\
&\leq \|\nabla(\Delta d-f(d))\|_2+(1+C\|\nabla d\|_2\|\Delta d\|_2)\|\nabla d\|_2\notag\\
&\leq \|\nabla(\Delta d-f(d))\|_2+C(1+\|\nabla d\|_2\|\Delta d-f(d)\|_2)\|\nabla d\|_2\notag.
\end{split}
\end{equation}
Therefore, applying H\"older inequality, Gagliardo-Nirenberg inequality and (\ref{esti:d-h3}) yields
\begin{equation}\label{esti:i20}
\begin{split}
|I_2|\leq& \epsilon\|\Delta u\|_2^2+C\|\nabla d\|_4^2\|\Delta d\|_4^2\\
\leq& \epsilon\|\Delta u\|_2^2+C\|\nabla d\|_2^{1/2}\|\Delta d\|_2^2\|\nabla\Delta d\|_2^{3/2}\\
\leq&\epsilon\|\Delta u\|_2^2+\epsilon\|\nabla\Delta d\|_2^2+C\|\nabla d\|_2^2\|\Delta d\|_2^8\\
\leq& \epsilon\|\Delta u\|_2^2+\epsilon\|\nabla(\Delta d-f(d))\|_2^2+ C\|\nabla d\|_2^2(\|\Delta d-f(d)\|_2^8+1)\\
&+C\|\nabla d\|_2^2(1+\|\nabla d\|_2^2\|\Delta d-f(d)\|_2^2).
\end{split}
\end{equation}

From the equation of $\sigma$ in (\ref{eq:g}) we have
\begin{align}\notag
I_3=&-\mu_1\int_Q\nabla_j(d_kd_pA_{kp}d_id_j)\nabla_l\nabla_lu_idx-\mu_2\int_Q\nabla_j(d_jN_i)\nabla_l\nabla_lu_idx\notag\\
&-\mu_3\int_Q\nabla_j(d_iN_j)\nabla_l\nabla_lu_idx
-\mu_5\int_Q\nabla_j(d_jd_kA_{ki})\nabla_l\nabla_lu_idx\notag\\
&-\mu_6\int_Q\nabla_j(d_id_kA_{kj})\nabla_l\nabla_lu_idx\notag.
\end{align}
Integration by parts yields
\begin{align}\notag%\label{eq:mu1}
&-\mu_1\int_Q\nabla_j(d_kd_pA_{kp}d_id_j)\nabla_l\nabla_lu_idx\\
=&\mu_1\int_Q(d_kd_pA_{kp}d_id_j)\nabla_l\nabla_l(A_{ij}+\Omega_{ij})dx\notag\\
=&\mu_1\int_Q(d_kd_pA_{kp}d_id_j)\nabla_l\nabla_lA_{ij}dx\notag
\end{align}
where we used the fact that $\Omega$ is antisymmetric which implies
\begin{align}\notag
\mu_1\int_Q(d_kd_pA_{kp}d_id_j)\nabla_l\nabla_l\Omega_{ij}dx
&=-\mu_1\int_Q\nabla_l(d_kd_pA_{kp}d_id_j)\nabla_l\Omega_{ij}dx\notag\\
&=\mu_1\int_Q\nabla_l(d_kd_pA_{kp}d_id_j)\nabla_l\Omega_{ji}dx\notag\\
&=\mu_1\int_Q\nabla_l(d_kd_pA_{kp}d_jd_i)\nabla_l\Omega_{ij}dx\notag\\
&=-\mu_1\int_Q(d_kd_pA_{kp}d_jd_i)\nabla_l\nabla_l\Omega_{ij}dx\notag
\end{align}
hence $\int_Q(d_kd_pA_{kp}d_id_j)\nabla_l\nabla_l\Omega_{ij}dx=0$. Therefore
\begin{equation}\label{eq:mu1}
\begin{split}
&-\mu_1\int_Q\nabla_j(d_kd_pA_{kp}d_id_j)\nabla_l\nabla_lu_idx\\
=&-\mu_1\int_Q|d^T\nabla Ad|^2dx-\mu_1\int_Q\nabla(d\otimes d)d\otimes dA\nabla Adx\\
=&-\mu_1\int_Q|d^T\nabla Ad|^2dx+I_{31}.
\end{split}
\end{equation}

Integration by parts also yields
\begin{equation}\label{eq:mu23}
\begin{split}
&-\mu_2\int_Q\nabla_j(d_jN_i)\nabla_l\nabla_lu_idx
-\mu_3\int_Q\nabla_j(d_iN_j)\nabla_l\nabla_lu_idx\\
=&\mu_2\int_Qd_jN_i\Delta(A_{ij}+\Omega_{ij})dx+\mu_3\int_Qd_iN_j\Delta(A_{ij}+\Omega_{ij})dx\\
=&(\mu_2+\mu_3)\int_Qd_jN_i\Delta A_{ij}dx+(\mu_2-\mu_3)\int_Qd_jN_i\Delta\Omega_{ij}dx\\
=&I_{32}+I_{33}
\end{split}
\end{equation} 
\begin{equation}\label{eq:mu56}
\begin{split}
&-\mu_5\int_Q\nabla_j(d_jd_kA_{ki})\nabla_l\nabla_lu_idx
-\mu_6\int_Q\nabla_j(d_id_kA_{kj})\nabla_l\nabla_lu_idx\\
=&(\mu_5+\mu_6)\int_Qd_jd_kA_{ki}\Delta A_{ij}dx+(\mu_5-\mu_6)\int_Qd_jd_kA_{ki}\Delta\Omega_{ij}dx\\
=&-2(\mu_5+\mu_6)\int_QAd\nabla d\nabla A\, dx
+(\mu_5-\mu_6)\int_QAd\Delta\Omega d\,dx\\
&-(\mu_5+\mu_6)\int_Q|d\nabla A|^2dx
=I_{34}+I_{35}+I_{36}.
\end{split}
\end{equation} 
Using integration by parts and equation (\ref{eqL3}), we infer
\begin{equation}\label{eq:i5}
\begin{split}
I_5=&\int_Q(\Delta d-f(d))\Delta \Omega d\, dx+2\int_Q(\Delta d-f(d))\nabla \Omega\nabla d\, dx\\
&+\int_Q(\Delta d-f(d))\Omega\Delta d\, dx\\
=&-\lambda_1(N,\Delta \Omega d)-\lambda_2(Ad,\Delta \Omega d)+2\int_Q(\Delta d-f(d))\nabla \Omega\nabla d\, dx\\
&-\int_Q(\Delta d-f(d))\nabla\Omega\nabla d\, dx-\int_Q\nabla(\Delta d-f(d))\Omega\nabla d\, dx\\
=&-\lambda_1(N,\Delta \Omega d)-\lambda_2(Ad,\Delta \Omega d)\\
&+\int_Q(\Delta d-f(d))\nabla\Omega\nabla d\, dx-\int_Q\nabla(\Delta d-f(d))\Omega\nabla d\, dx\\
=&I_{51}+I_{52}+I_{53}+I_{54}
\end{split}
\end{equation}
\begin{equation}\label{eq:i6}
\begin{split}
I_6=&\lambda_2\int_QN\Delta(Ad)dx+\frac{\lambda_2^2}{\lambda_1}\int_QAd\Delta(Ad)dx\\
=&\lambda_2\int_QN\Delta Ad\,dx+2\lambda_2\int_QN\nabla A\nabla d\,dx\\
&+\lambda_2\int_QNA\Delta d\,dx-\frac{\lambda_2^2}{\lambda_1}\int_Q|d\nabla A|^2dx
-\frac{\lambda_2^2}{\lambda_1}\int_Q|A\nabla d|^2dx\\
=&I_{61}+I_{62}+I_{63}+I_{64}+I_{65}.
\end{split}
\end{equation}
Note that there are cancelations among $I_3, I_5$ and $I_6$, due to the parameter relations (\ref{re:con}) and (\ref{eq:Par}), and the assumption (\ref{ass:para2}). Indeed, 
$$
I_{32}+I_{61}=0, \, \, I_{33}+I_{51}=0, \, \, I_{35}+I_{52}=0, \, \, I_{36}+I_{64}\leq 0.
$$
In the following we estimate the rest terms in $I_3, I_5$ and $I_6$. Applying H\"older's inequality, Gagliardo-Nirenberg's inequality and Agmon's inequality (\ref{Ag}) yields
\begin{equation}\label{esti:i31}
\begin{split}
|I_{31}|\leq&\epsilon\|\Delta u\|_2^2+C\|d\|_\infty^6\|\nabla u\|_4^2\|\nabla d\|_4^2\\
\leq&\epsilon\|\Delta u\|_2^2+C\|\nabla d\|_2^3\|\Delta d\|_2^3\|\nabla u\|_2^{1/2}\|\Delta u\|_2^{3/2}
\|\nabla d\|_2^{1/2}\|\Delta d\|_2^{3/2}\\
\leq&\epsilon\|\Delta u\|_2^2+C\|\nabla d\|_2^{14}\|\Delta d\|_2^{18}\|\nabla u\|_2^{2}\\
\leq&\epsilon\|\Delta u\|_2^2+C\|\nabla d\|_2^{14}(\|\Delta d-f(d)\|_2^{18}+1)\|\nabla u\|_2^{2}
\end{split}
\end{equation}
\begin{equation}\label{esti:i34}
|I_{34}|
\leq\epsilon\|\Delta u\|_2^2+C\|\nabla d\|_2^{6}(\|\Delta d-f(d)\|_2^{10}+1)\|\nabla u\|_2^{2}
\end{equation}
\begin{equation}\label{esti:i53}
|I_{53}|
\leq\epsilon(\|\Delta u\|_2^2+\|\nabla(\Delta d-f)\|_2^2)+C\|\nabla d\|_2^{2}(\|\Delta d-f(d)\|_2^{6}+1)\|\Delta d-f(d)\|_2^{2}
\end{equation}
\begin{equation}\label{esti:i54}
|I_{54}|
\leq\epsilon(\|\Delta u\|_2^2+\|\nabla(\Delta d-f)\|_2^2)+C\|\nabla d\|_2^{2}(\|\Delta d-f\|_2^{6}+1)\|\nabla u\|_2^{2}
\end{equation}
and 
\begin{equation}\label{esti:i62}
\begin{split}
|I_{62}|\leq&\int_Q|Ad\Delta u\nabla d|dx+\int_Q|(\Delta d-f(d))\Delta u\nabla d|dx\\
\leq&\epsilon(\|\Delta u\|_2^2+\|\nabla(\Delta d-f)\|_2^2)+C\|\nabla d\|_2^{2}\\
&\cdot(\|\Delta d-f\|_2^{10}+\|\Delta d-f\|_2^{6}+1)\left(\|u\|_2^2+\|\Delta d-f\|_2^{2}\right)
\end{split}
\end{equation}
since the two integrals are similar to $I_{34}$ and $I_{53}$ respectively. While, we have by integration by parts
\begin{equation}\notag%\label{esti:i63}
\begin{split}
I_{63}+I_{65}=&-\lambda_2\int_QN\nabla A\nabla d\,dx-\lambda_2\int_Q\nabla NA\nabla d\, dx-\frac{\lambda_2^2}{\lambda_1}\int_Q|A\nabla d|^2\, dx\\
=&-\lambda_2\int_QN\nabla A\nabla d\,dx+\frac{\lambda_2^2}{\lambda_1}\int_Q\nabla (Ad)A\nabla d\, dx\\
&+\frac{\lambda_2}{\lambda_1}\int_Q\nabla (\Delta d-f)A\nabla d\, dx-\frac{\lambda_2^2}{\lambda_1}\int_Q|A\nabla d|^2\, dx\\
=&-\lambda_2\int_QN\nabla A\nabla d\,dx+\frac{\lambda_2^2}{\lambda_1}\int_Q\nabla AdA\nabla d\, dx\\
&+\frac{\lambda_2}{\lambda_1}\int_Q\nabla (\Delta d-f)A\nabla d\, dx
\end{split}
\end{equation}
where the three integrals are similar to $I_{62}, I_{34}$ and $I_{54}$ respectively. Thus
\begin{equation}\label{esti:i63}
\begin{split}
|I_{63}+I_{65}|\leq&\epsilon(\|\Delta u\|_2^2+\|\nabla(\Delta d-f)\|_2^2)+C\|\nabla d\|_2^{2}\\
&\cdot(\|\Delta d-f\|_2^{10}+\|\Delta d-f\|_2^{6}+1)\left(\|u\|_2^2+\|\Delta d-f\|_2^{2}\right).
\end{split}
\end{equation}

The estimate for $I_4+I_8$ is as follows. The facts $\nabla d\cdot f(d)=\nabla F(d)$ and $\nabla\cdot u=0$ imply that $(f,\Delta u\cdot\nabla d)=0$. Also by $\nabla\cdot u=0$, it has $(\Delta d-f, u\cdot\nabla(\Delta d-f))=0$. Thus, we have
\begin{equation}\notag
\begin{split}
I_4&=-(\Delta d-f,\Delta u\cdot\nabla d)-2\int_\Omega(\Delta d_i-f_i)\nabla_lv_j\nabla_l\nabla_jd_idx-(\Delta d-f,u\cdot\nabla\Delta d)\\
&=-(\Delta d\nabla d,\Delta u)+2\int_\Omega\nabla_j(\Delta d_i-f_i)\nabla_lv_j\nabla_ld_idx-(\Delta d-f,u\cdot\nabla f(d))\\
&=-(\Delta d\nabla d,\Delta u)+2\int_\Omega\nabla_j(\Delta d_i-f_i)\nabla_lv_j\nabla_ld_idx-I_8.
\end{split}
\end{equation}
Therefore, 
\begin{equation}\notag
I_4+I_8=-(\Delta d\nabla d,\Delta u)+2\int_\Omega\nabla_j(\Delta d_i-f_i)\nabla_lv_j\nabla_ld_idx.
\end{equation}
where the two terms are similar to $I_2, I_{54}$ respectively.
Hence, from (\ref{esti:i20}) and (\ref{esti:i54}) we have
\begin{equation}\label{esti:i48}
|I_4+I_8|\leq\epsilon\left(\|\Delta u\|_2^2+\|\nabla(\Delta d-f)\|_2^2\right)
+C\|\nabla d\|_2^2\left(\Phi^8+\Phi^6+1\right).
\end{equation}

By Agmon's inequality, the terms $I_7$ and $I_9$ are estimated as
\begin{equation}\label{esti:i7}
\begin{split}
|I_7|&\leq C\|f'(d)\|_\infty\|\Delta d-f(d)\|_2^2\leq C\|d\|_\infty^2\|\Delta d-f(d)\|_2^2\\
&\leq C\|\nabla d\|_2(\|\Delta d-f(d)\|_2+1)\|\Delta d-f(d)\|_2^2,
\end{split}
\end{equation}
\begin{equation}\label{esti:i9}
\begin{split}
|I_9|&\leq C\|f'(d)\|_\infty\|d\|_\infty\left(\|\Delta d-f(d)\|_2^2+\|\nabla u\|_2^2\right)\\
&\leq C\|\nabla d\|_2^{3/2}(\|\Delta d-f(d)\|_2^{3/2}+1)\Phi^2.
\end{split}
\end{equation}
Denote $\mathcal D=\|\Delta d-f(d)\|_2$.  Combining (\ref{eq:phi1}), (\ref{esti:i1}) and (\ref{esti:i20})-(\ref{esti:i9}) yields
\begin{equation}\label{eq:phi2}
\begin{split}
&\frac{1}{2}\frac{d}{dt}\Phi^2(t)+\left(\mu_4-\epsilon\right)\int_Q|\Delta u|^2dx\\
&+\left(-\frac{1}{\lambda_1}-\epsilon\right)\int_\Omega|\nabla(\Delta d-f(d))|^2dx+
\mu_1\int_Q|d^T\nabla Ad|^2dx\\
\leq&CM_1^{-4}\|u\|_2^2\|\sqrt{\rho}\nabla u\|_2^8+C\|\nabla d\|_2^2(\mathcal D^8+\mathcal D^2+1)\\
&+C\|\nabla d\|_2^2(\mathcal D^{18}+\mathcal D^{10}+\mathcal D^6+1)(\|\nabla u\|_2^2+\mathcal D^2)\\
&+C\|\nabla d\|_2(\mathcal D+1)\mathcal D^2+C\|\nabla d\|_2^{3/2}(\mathcal D^{3/2}+1)(\|\nabla u\|_2^2+\mathcal D^2)\\
%\leq&CM_1^{-4}\|u\|_2^2\|\sqrt{\rho}\nabla u\|_2^8+C\|\nabla d\|_2(\mathcal D^8+\mathcal D^2+1)\notag\\
%&+C\|\nabla d\|_2^2(\mathcal D^{18}+\mathcal D^{10}+\mathcal D^6+1)(M_1^{-1}\|\sqrt{\rho}\nabla u\|_2^2+\mathcal D^2)\notag\\
%&+C\|\nabla d\|_2(\mathcal D^2+\mathcal D+1)(M_1^{-1}\|\sqrt{\rho}\nabla u\|_2^2+\mathcal D^2)\notag\\
\leq&C\|\nabla d\|_2^2+(\|u\|_2+\|\nabla d\|_2)(\Phi^{18}+\Phi^{10}+\Phi^{6}+\Phi^{2}+\Phi+1)\Phi^{2},
\end{split}
\end{equation}
where we used the facts $\|\nabla u\|_2^2\leq M_1^{-1}\|\sqrt{\rho}\nabla u\|_2^2$ and $\|u\|_2+\|\nabla d\|_2\leq C$ by the basic energy estimate.

Set
$$
\tilde\Phi^2 = \Phi^2+\|u\|_2 + \|\nabla d\|_2
$$
and observe from (\ref{eq:phi2}) that
\begin{equation}\label{eq:phi3}
\frac d{dt} \tilde\Phi^2 \leq C(1+(\|u\|_2 + \|\nabla d\|_2)\tilde\Phi^{18})\tilde \Phi^2.
\end{equation}
Recall that $\|u\|_2 + \|\nabla d\|_2$ is small by the basic energy estimate (\ref{energy}) and the smallness assumption (\ref{ass:small}). Thus following a similar argument as in \cite{DQSr} we infer that 
for a small enough $\epsilon_0>0$ in the  assumption (\ref{ass:small}), the conclusions of the Lemma will be proved. 
Suppose that
\begin{equation}\label{sup:small}
(\frac{M_2}{\mu_4}+1)\rho_0\|u_0\|_{H^1} ^2+ \|d_0\|_{H^1}^2 + \|\Delta d_0-f(d_0)\|_{2}^2= \epsilon_0.
\end{equation}
By the basic energy estimate we have
$$
\|u\|_{2} + \|\nabla d\|_{2}\leq \|u_0\|_{2} + \|\nabla d_0\|_{2}\leq\sqrt{2\epsilon_0}.
$$
We claim that, if $\epsilon_0$ is so small that
\begin{equation}\label{small}
\sqrt{2\epsilon_0} (4e^{2C}\sqrt{\epsilon_0})^{9} < 1
\end{equation}
then
$$
\tilde\Phi^2 < 4 e^{2C}\sqrt{\epsilon_0}, \, \, \mbox { for all } t>0.
$$
First we prove the claim for $t\in [0, 1]$. Assume otherwise, there must be $t_0\in (0, 1)$ such that
\begin{equation}\label{wrong}
\begin{cases}
\tilde\Phi^2(t_0) = 4 e^{2C}\sqrt{\epsilon_0}\\
\tilde\Phi^2(t) \leq 4 e^{2C}\sqrt{\epsilon_0}, \, \, \mbox { for all } t\in (0, t_0].
\end{cases}
\end{equation}
Therefore, from (\ref{eq:phi3}), by the choice of $\epsilon_0$ in \eqref{small}, we have
$$
\frac d{dt} \tilde\Phi^2 \leq 2C\tilde\Phi^2
$$
for all $t\in (0, t_0)$ and $\tilde\Phi^2(0) \leq \epsilon_0+\sqrt{2\epsilon_0}<3\sqrt{\epsilon_0}$, which implies that
$$
\tilde\Phi^2(t_0) \leq e^{2C}\tilde\Phi^2(0)\leq 3e^{2C}\sqrt{\epsilon_0}
$$
and thus contradicts \eqref{wrong}. For $t >1$, we simply observe, %as in \cite{DQSr},
 that from 
the basic energy inequality \eqref{energy} and the assumption (\ref{sup:small})
$$
\int_{t-1}^t \tilde\Phi^2 (t)dt \leq 3\sqrt{\epsilon_0}
$$
which implies that there is $t_0\in (t-1, t)$ such that
$$
\tilde\Phi^2(t_0) \leq 3\sqrt{\epsilon_0}.
$$
One may repeat the above argument to conclude that
$$
\tilde\Phi^2(t)\leq 4 e^{2C}\sqrt{\epsilon_0}, \, \, \mbox { for all } t>0.
$$
The inequality (\ref{L2H2}) follows directly from the basic energy estimate (\ref{energy}) and (\ref{eq:phi3}). It completes the proof of the lemma.
\cbdu

Without the assumption of small initial data, the higher order energy estimate holds for short time.

\begin{Lemma}\label{le:phi-short} (3D short time.) Let $Q\subset\mathbb{R}^3$ and   the initial data $(\rho_0,u_0,d_0)$  satisfy (\ref{initden})-(\ref{initu}). Suppose  that $\Phi(0)^2=\|\sqrt{\rho_0}\nabla u_0\|_{2}^2+\|\Delta d_0-f(d_0)\|_{2}^2<\infty$. Then the approximating solutions $(\rho^m,u^m,d^m)$ obtained in Lemma \ref{le:weak}  satisfy (\ref{phiest}) and (\ref{L2H2}) on the time interval $(0,T]$ for some $T<\delta$.
\end{Lemma}
\pf
Set
$$
\tilde\Phi_m^2 = \Phi_m^2+\|u^m\|_{L^2} + \|\nabla d^m\|_{L^2}+1.
$$
By (\ref{eq:phi2}) and the uniform estimate $\|u^m\|_2+\|\nabla d^m\|_2\leq C$, we have
\begin{equation}\label{eq:phi3}
\frac d{dt} \tilde\Phi_m^2 \leq C \tilde\Phi_m^{20}.
\end{equation}
Hence,
\begin{equation}\notag
\int\frac {d\tilde\Phi_m^2}{\tilde\Phi_m^{18}} \leq C \int_0^T\tilde\Phi_m^{2}dt
\end{equation}
which produces for any $t\in(0,T]$
\begin{equation}\notag
\tilde\Phi_m^{16}(t)\leq\frac{\tilde\Phi_m^{16}(0)}{1-C\tilde\Phi_m^{16}(0)\int_0^T\tilde\Phi_m^{2}(0)dt}.
\end{equation}
The basic energy estimate (\ref{energy}) indicates, there exists a small time $T<\delta$ such that 
$$
C\tilde\Phi_m^{16}(0)\int_0^T\tilde\Phi_m^{2}(0)dt\leq\frac 12
$$
which consequently implies 
$$
\tilde\Phi_m^2(t)\leq C\tilde\Phi_m^2(0), \, \, \mbox { for all } t\in(0,T].
$$
The conclusion of the lemma thus holds true.
\cbdu

Next we shall establish the higher order energy estimate with general data in two dimensional case. Note that in 2D, we usually assume in addition either
$$
(i)\, \, \mu_1\geq 0, \, \, \lambda_2=0, \, \, \mbox { or } \, \, (ii)\, \, \mu_1=0, \, \, \lambda_2\neq 0.
$$
For the literature of the Ericksen-Leslie system with assumption (i), we refer the readers to \cite{LL1}; while for the work with assumption (ii), we refer to \cite{CR, GW, SL, WXL}.

\begin{Lemma}\label{le:phi2d} (2D.) Let $Q\subset\mathbb{R}^2$ and the initial data $(\rho_0,u_0,d_0)$  satisfy (\ref{initden})-(\ref{initu}). Suppose  that $\Phi(0)^2=\|\sqrt{\rho_0}\nabla u_0\|_{2}^2+\|\Delta d_0-f(d_0)\|_{2}^2<\infty$. Assume either (i) or (ii) holds. Then the approximating solutions $(\rho^m,u^m,d^m)$ obtained in Lemma \ref{le:weak}  satisfy 
\begin{equation}\label{phiest-2d}
\|\sqrt{\rho^m}\nabla u^m\|_{L^2}^2 + \|\Delta d^m-f(d^m)\|_{L^2}^2 \leq (\Phi^2(0)+1)e^{C(\rho_0\|u_0\|_2^2+\|\nabla d_0\|_2^2+1)}
\end{equation}
with $t\in [0, T]$ and
\bg\label{L2H2-2d}
\int_0^T\|\Delta u^m\|_{2}^2 +\|\nabla(\Delta d^m-f(d^m))\|_{2}^2+\mu_1\|d^{mT}\nabla A^md^m\|_2^2dt \leq C
\ed
for all $T>0$.
The constants $C$ depend only on initial data, $Q$, $M_1$, $M_2$ and the physical coefficients in system (\ref{eqL1})-(\ref{eqL4}).
\end{Lemma}
\pf
Case (i). Since $\lambda_2=0$, a maximum principle for $d$ holds (cf. \cite{LL1}).  The Agmon's inequality is not needed here. The highest nonlinear term $I_{31}$ is estimated as
\begin{align}\notag
|I_{31}|\leq &C\|d^m\|_\infty^3\int_Q|\nabla d^m\nabla u^m\Delta u^m|dx\\
\leq &\epsilon\|\Delta u^m\|_2^2+C\|\nabla u^m\|_4^2\|\nabla d^m\|_4^2\notag\\
\leq &\epsilon\|\Delta u^m\|_2^2+C\|\nabla u^m\|_2^2(\|\Delta d^m-f(d^m)\|_2^2+1)\notag
\end{align}
where we used the 2D Gagliardo-Nirenberg inequality in Lemma \ref{le:GN}. Similarly, we apply H\"older's equality, Gagliardo-Nirenberg inequality, the basic energy estimate and 
$\|d^m\|_\infty\leq C$ to all the other terms in (\ref{eq:phi0}) and obtain that
\begin{equation}%\label{eq:phi-2d1}
\begin{split}
&\frac{1}{2}\frac{d}{dt}\Phi_m^2(t)+\left(\mu_4-\epsilon\right)\int_Q|\Delta u^m|^2dx\\
&+\left(-\frac{1}{\lambda_1}-\epsilon\right)\int_\Omega|\nabla(\Delta d^m-f(d^m))|^2dx+
\mu_1\int_Q|d^{mT}\nabla A^md^m|^2dx\\
\leq&C(\Phi_m^{4}+\Phi_m^{2}+1).
\end{split}
\end{equation}
Let $\tilde\Phi_m^2=\Phi_m^2+1$. We have 
\begin{equation}\notag%\label{eq:phi5}
\begin{split}
\frac{d}{dt}\tilde\Phi_m^2(t)+\mu_4\int_Q|\Delta u^m|^2dx
-\frac{1}{\lambda_1}\int_\Omega|\nabla(\Delta d^m-f(d^m))|^2dx\\
+2\mu_1\int_Q|d^{mT}\nabla A^md^m|^2dx
\leq C\tilde\Phi_m^{4}.
\end{split}
\end{equation}
Thus, 
$$
\frac{1}{\tilde\Phi_m^2}\frac d{dt} \tilde\Phi_m^2 \leq C\tilde\Phi_m^2.
$$
The conclusion of the lemma follows from the above inequality and the basic energy estimate
(\ref{energy}) immediately. \\

Case (ii). Without loss of the generality, we assume $|d^m|$ is large. Since otherwise we have $\|d^m\|_\infty\leq C$ and hence it can be handled similarly as in case (i). Thus there exists a constant $C_1$ such that 
\begin{equation}\notag
\int_Q|d^m|^2dx\leq |Q|^{1/2}\left(\int_Q|d^m|^4dx\right)^{1/2}\leq C_1|Q|^{1/2}\left(\int_Q(|d^m|^2-1)^2dx\right)^{1/2}\leq C
\end{equation}
where we used the fact $\int_QF(d^m)dx$ is bounded by the initial data from the basic energy estimate.
It then follows from the Agmon's inequality in 2D that
\begin{equation}\notag
\|d^m\|_\infty^2\leq C\|d^m\|_2\|\Delta d^m\|_2\leq C\|\Delta d^m\|_2.
\end{equation}
Hence, slightly different from (\ref{esti:d-h3}) we have
\begin{equation}\notag
\|\nabla\Delta d^m\|_2
\leq \|\nabla(\Delta d^m-f(d^m))\|_2+C(1+\|\Delta d^m-f(d^m)\|_2)\|\nabla d^m\|_2.
\end{equation}
Keep in mind that $\mu_1=0$ and thus the highest nonlinear term disappears here, and we have the Gagliardo-Nirenberg inequality in 2D. An analogous computation as in the proof of Lemma \ref{le:phi} produces that
\begin{equation}\notag%\label{eq:phi4}
\begin{split}
&\frac{1}{2}\frac{d}{dt}\Phi_m^2(t)+\left(\mu_4-\epsilon\right)\int_Q|\Delta u^m|^2dx\\
&+\left(-\frac{1}{\lambda_1}-\epsilon\right)\int_\Omega|\nabla(\Delta d^m-f(d^m))|^2dx\\
\leq&\epsilon\|\sqrt{\rho^m}\nabla u^m\|_2^2\|\Delta u^m\|_2^2+C(\Phi_m^{4}+\Phi_m^{2}+1).
\end{split}
\end{equation}
Set $\tilde\Phi_m^2=\Phi_m^2+1$. It follows 
\begin{equation}\label{eq:phi5}
\frac{d}{dt}\tilde\Phi_m^2(t)\leq-\left(\mu_4-\epsilon\tilde\Phi_m^2\right)\int_Q|\Delta u^m|^2dx+C\tilde\Phi_m^{4}.
\end{equation}
For a mall enough $\epsilon$ such that 
\begin{equation}\label{small-e}
\epsilon 2\tilde\Phi^2(0)e^{C(\rho_0\|u_0\|_2^2+\|\nabla d_0\|_2^2+1)}\leq \frac{\mu_4}{2}
\end{equation}
then
$$
\tilde\Phi_m^2(t) < 2 \tilde\Phi^2(0)e^{C(\rho_0\|u_0\|_2^2+\|\nabla d_0\|_2^2+1)}, \, \, \mbox { for all } t>0.
$$
First we prove the claim for $t\in [0, 1]$. Assume otherwise, there must be $t_0\in (0, 1)$ such that
\begin{equation}\label{wrong-2d}
\begin{cases}
\tilde\Phi_m^2(t_0) = 2 \tilde\Phi^2(0)e^{C(\rho_0\|u_0\|_2^2+\|\nabla d_0\|_2^2+1)}\\
\tilde\Phi_m^2(t) \leq 2 \tilde\Phi^2(0)e^{C(\rho_0\|u_0\|_2^2+\|\nabla d_0\|_2^2+1)}, \, \, \mbox { for all } t\in (0, t_0].
\end{cases}
\end{equation}
Therefore, from (\ref{eq:phi5}), by the choice of $\epsilon$ in \eqref{small-e}, we have
$$
\frac{1}{\tilde\Phi_m^2}\frac d{dt} \tilde\Phi_m^2 \leq C\tilde\Phi_m^2
$$
for all $t\in (0, t_0)$ %and $\tilde\Phi^2(0) \leq \epsilon_0+\sqrt{2\epsilon_0}<3\sqrt{\epsilon_0}$, 
which implies that
$$
\tilde\Phi_m^2(t_0) \leq \tilde\Phi^2(0)e^{C\int_0^1\tilde\Phi_m^2(t)dt}\leq \tilde\Phi^2(0)e^{C(\rho_0\|u_0\|_2^2+\|\nabla d_0\|_2^2+1)}
$$
and thus contradicts \eqref{wrong-2d}. For $t >1$, we apply a similar argument as in the proof of Lemma \ref{le:phi} to show that
$$
\tilde\Phi_m^2(t)\leq 2 \tilde\Phi^2(0)e^{C(\rho_0\|u_0\|_2^2+\|\nabla d_0\|_2^2+1)}. %\, \, \mbox { for all } t>0.
$$
%The inequality (\ref{L2H2}) follows directly from the basic energy estimate (\ref{energy}) and (\ref{eq:phi3}). 
It completes the proof of the lemma.
\cbdu

\subsection{Passage of limit to the weak solutions.} 
\label{sec:weak-sol}

In Subsection \ref{sec:higher-energy}, the higher order energy estimates in Lemmas \ref{le:phi}, \ref{le:phi-short} and \ref{le:phi2d} are all independent of the approximation index $m$ and time $t$. It implies the following uniform estimates hold
$$
u^m\in L^\infty(0,T; V)\cap L^2(0,T; H_p^2), 
$$
$$
d^m\in L^\infty(0,T; H_p^2)\cap L^2(0,T; H_p^3).
$$
By the mass conservation equation (\ref{eqL1}), we also have
$$
0<M_1\leq\rho^m\leq M_2.
$$
Moreover, the positive lower bound of density $\rho^m$ combined with the uniform estimates for $u^m$ and $d^m$ indicates that
$$
u_t^m\in L^2(0,T; L_p^2), \, \, \,
d_t^m\in L^2(0,T; H_p^1).
$$
Therefore, with these higher order energy estimates, a standard procedure will show that the approximating solutions $(\rho^m, u^m, d^m)$ (a subsequence) converge to a limit $(\rho, u, d)$
such that the limit is a weak solution of system (\ref{eqL1})-(\ref{eqL4}) and satisfies 
\bg\label{esti:weak-rho}
0<M_1\leq\rho\leq M_2,
\ed
\bg\label{esti:weak-u}
u\in L^\infty(0,T; V)\cap L^2(0,T; H_p^2), 
\ed
\bg\label{esti:weak-d}
d\in L^\infty(0,T; H_p^2)\cap L^2(0,T; H_p^3).
\ed
Thus we obtained the existence of weak solutions.
\begin{Theorem}(3D)\label{Thm: weak-3D}
Let $Q\subset\mathbb{R}^3$ and the initial data $(\rho_0,u_0,d_0)$  satisfy the conditions (\ref{initden})-(\ref{initu}). Suppose  that $\Phi^2(0)=\|\sqrt{\rho_0}\nabla u_0\|_{L^2}^2+\|\Delta d_0-f(d_0)\|_{L^2}^2<\infty$. \\
(I) There is $\epsilon_0>0$ such that if
\begin{equation}\notag%\label{ass:small}
(\frac{M_2}{\mu_4}+1)\rho_0\|u_0\|_{H^1}^2 + \|d_0\|_{H^1}^2+\|\Delta d_0-f(d_0)\|_2^2 \leq \epsilon_0,
\end{equation}
then system (\ref{eqL1})-(\ref{eqL4}) with the periodic boundary condition (\ref{bd}) has a weak  solution $(\rho, u, d)$ satisfying the basic energy estimate inequality 
%{\color{red}(change it to integral form, no derivative on time)}
\begin{align}\label{energy-weak}
\mathcal E(t)&+\int_{t_0}^t\int_Q\mu_1|(d)^TAd|^2+\mu_4|\nabla u|^2-\frac{1}{\lambda_1}|\Delta d-f(d)|^2dx\\
&+\left(\mu_5+\mu_6+\frac{\lambda_2^2}{\lambda_1}\right)\int_{t_0}^t\int_Q|Ad|^2dx\leq \mathcal E(t_0),\notag
\end{align}
for all $t>t_0$ and almost every $t_0$. Moreover, the weak solution $(\rho, u, d)$ satisfies the higher order energy estimates (\ref{esti:weak-rho})-(\ref{esti:weak-d}) for all $T>0$. \\
(II) On the other hand, without the smallness assumption on initial data, there exists a weak solution to the system satisfying (\ref{energy-weak}) for all $t>t_0$ and a.e. $t_0$, but satisfying (\ref{esti:weak-rho})-(\ref{esti:weak-d}) only for a short time $T<\delta$, for some $\delta>0$.
%\bg\notag
%\mathcal E^m(t)=\int_Q\frac{1}{2}\rho^m |u^m|^2+\frac{1}{2}|\nabla d^m|^2+\rho^m \mathcal F(d^m)dx.
%\ed
\end{Theorem}

\begin{Theorem}(2D)\label{Thm: weak-2D}
Let $Q\subset\mathbb{R}^2$ and the initial data $(\rho_0,u_0,d_0)$  satisfy the conditions (\ref{initden})-(\ref{initu}). Suppose  that $\Phi^2(0)=\|\sqrt{\rho_0}\nabla u_0\|_{L^2}^2+\|\Delta d_0-f(d_0)\|_{L^2}^2<\infty$. Assume in addition that either
$$
(i)\, \, \mu_1\geq 0, \, \, \lambda_2=0, \, \, \mbox { or } \, \, (ii)\, \, \mu_1=0, \, \, \lambda_2\neq 0.
$$
The system (\ref{eqL1})-(\ref{eqL4}) with the periodic boundary condition (\ref{bd}) has a weak  solution $(\rho, u, d)$ satisfying the basic energy estimate inequality
(\ref{energy-weak}) with $t\in(0,T]$ and the higher order energy estimates (\ref{esti:weak-rho})-(\ref{esti:weak-d}) for all the time $T>0$.
\end{Theorem}

Applying a similar method to obtain the higher order energy estimates as in Subsection \ref{sec:higher-energy}, the weak solution $(\rho, u, d)$ satisfies even higher order energy estimates provided more regular initial data.

\begin{Corollary}
Let $Q\subset \R^n$ with $n=2,3$. The initial data $(\rho_0, u_0, d_0)$ satisfies the conditions in Theorem \ref{Thm: weak-3D} and Theorem \ref{Thm: weak-2D}, respectively for $n=3, 2$.  Assume, in addition, $u_0\in H_p^2$ and $d_0\in H_p^3$. Then there exists a weak solution $(\rho, u, d)$ to system (\ref{eqL1})-(\ref{eqL4}) satisfying the basic energy estimate (\ref{energy-weak}) and (\ref{esti:weak-rho}). Mean while,
\bg\label{esti:weak-u-h}
u\in L^\infty(0,T; H_p^2)\cap L^2(0,T; H_p^3), 
\ed
\bg\label{esti:weak-d-h}
d\in L^\infty(0,T; H_p^3)\cap L^2(0,T; H_p^4).
\ed
\end{Corollary}

Indeed, compared to the simplified Ericksen-Leslie model studied in \cite{DQSr}, the full model contains higher order derivative terms, for example, the term $\nabla\cdot(Nd)$ in equation (\ref{eqL2}). When applying the $L^p$ theory of parabolic equation on (\ref{eqL2}) to improve the regularity of velocity, the higher order energy estimates (\ref{esti:weak-u-h}) and (\ref{esti:weak-d-h}) are needed. Please see a detailed discussion in Subsection \ref{sec:classical}.
% {\color{red}double check this argument.}

\subsection{Auxiliary estimates on density}
\label{sec:density}
As a consequence of the higher order energy estimates obtained in Subsection \ref{sec:weak-sol} for the weak solution $(\rho, u, d)$, the regularity of the density can be improved in both of two and three dimensional cases.

\begin{Lemma}\label{le:density2D}  (\cite{AKM}) 
Assume that $\rho(0)\in C^1(Q)$ and that $Q\subset\R^2$ is smooth and bounded.
 Suppose $u\in L^\infty(0, T; H^1)\bigcap L^2(0, T; H^2)$, and 
$$
\rho_t + u\cdot\nabla\rho = 0
$$
in $Q\times (0, T)$. Then $\rho\in C^\alpha(\Omega\times[0,T])$ for some $\alpha\in (0,1)$ which depends only on the initial data, $T$ and $Q$.
\end{Lemma}

%The  proof of this lemma in \cite{AKM} used a sharp embedding inequality (cf.  \cite{Il})
%\bg\notag
%|f(x_1)-f(x_2)| \leq  C\|f\|_{H^2(Q)}|x_1-x_2|(1+|\ln|x_1-x_2||),
%\ed
%where $x_1, x_2\in Q$ and $C$ is a constant depending on $Q$.

\begin{Lemma}\label{le:density3D}  (\cite{DQSr}) 
Assume that $\rho(0)\in C^1(Q)$ and that $Q\subset\R^3$ is smooth and bounded.
 Suppose $u\in L^\infty(0, T; H^1)\bigcap L^2(0, T; H^2)$, and 
$$
\rho_t + u\cdot\nabla\rho = 0
$$
in $Q\times (0, T)$. Let $t_1\in(0,T)$ and $p\in Q$, define
$$
\Aa_{(p,t_1)}=(B_{r_0}(p)\cap Q)\times([t_1-r_0,t_1+r_0]\cap[0,T]).
$$
Then, for any $\epsilon>0$, there exists  $r_0>0$ such that for $p\in Q$ and all $T>0$,
\bg\notag
\sup_{(q,t_2)\in\Aa_{(p,t_1)}}|\rho(q, t_2)- \rho(p, t_1)| \leq \epsilon.
\ed
\end{Lemma}
We refer the readers to \cite{DQSr} for a detailed proof of this lemma.

\begin{Remark}
In 2D, the H\"{o}lder continuity for the fluid density guarantees that a frozen coefficient method (cf. \cite{LS1}) can be applied to the Navier-Stokes equation (\ref{eqL2}) and hence the regularity of the velocity $u$ will be improved through a standard $L^p$ theory for the parabolic equation.  In 3D, the density has the property of small oscillations over small balls in $Q\times[0,T]$ provided that either the initial data is small or for short time. This turns out to be enough to carry out the frozen coefficient method to improve the regularity of the fluid velocity too. We refer the reader to \cite{LS1} for a general idea of the frozen coefficient method. Also a more specific discussion relevant to the current model can be found in \cite{DQSr} (Appendix 6).
\end{Remark}

\subsection{Classical solution}
\label{sec:classical}

%\subsection{Proof of Theorem \ref{Mthm2}}

In this subsection, using the the estimates (\ref{esti:weak-rho})-(\ref{esti:weak-d}),  we apply the frozen coefficient method to improve regularities for $\rho, p, u, d$ by a  bootstrapping argument among the three equations (\ref{eqL1})-(\ref{eqL3}). Hence Theorem \ref{Mthm1} and \ref{Mthm2} will be proved. The process of obtaining regular solutions in two dimensional case is much easier than the one in three dimensional case. And based on the previous work for the simplified Ericksen-Leslie model in \cite{DQSr}, we just briefly show the steps to obtain regular solutions in 3D in the following.  \\

We notice from (\ref{esti:weak-u}) and (\ref{esti:weak-d})
$$
u\in L^\infty(0,T;L_p^6), \ \ \nabla d\in L^\infty(0,T;L_p^6),
$$
which implies
$$
u\cdot\nabla d\in L^\infty(0,T;L_p^3).
$$
In the mean time, (\ref{esti:weak-u}) and (\ref{esti:weak-d}) indicate
$$
\Omega d, Ad\in L^\infty(0,T;L_p^q), \, \, \mbox { for all } q>1.
$$
By the standard parabolic estimates on equation (\ref{eqL3}) (cf. \cite{LS} and \cite{Am}), we have
$$
d\in W^{1,r}(W_p^{2,3}), \, \, \mbox { for all } r>1
$$
which implies that $\nabla d\in L^\infty(0,T;L_p^q)$ for any $q\in(1,\infty)$. Thus, we have
$$
u\cdot\nabla d\in L^\infty(0,T;L_p^q), \ \ \forall q\in(1,6).
$$
Applying the same standard parabolic estimates on (\ref{eqL3}) again yields
\bg\label{dw2r}
d\in W^{1, r}(W_p^{2,q}), \forall r\in (1, \infty) \ \text{and} \ q \in(1,6),
\ed
which implies that $d\in C^{\alpha/2,1+\alpha}([0,T]\times\bar Q)$ for some $\alpha \in (0, 1)$ and
$$
\nabla d\Delta d \in L^\infty(0, T; L_p^q), \quad \forall q \in (1, 6).
$$

In Navier-Stokes equation (\ref{eqL2}), the estimates for the conservation of momentum with constant density can be extended to the non-constant density case when Lemma \ref{le:density3D} is available. This is done via the frozen coefficient method.

The estimates (\ref{esti:weak-u}) and (\ref{esti:weak-d}) implies
$$
u\cdot\nabla u\in L^\infty(0,T;L_p^{\frac{3}{2}}),
$$
while (\ref{esti:weak-u-h}) and (\ref{esti:weak-d-h}) implies
$$
\nabla\cdot\sigma\in L^\infty(0,T;L_p^{\frac{3}{2}}).
$$
Now we apply the frozen coefficient method using the oscillation estimates for the density as in Lemma \ref{le:density3D}  to yield
$$
u\in W^{1, q}(W_p^{2,3/2}), \quad \forall q\in (1, \infty).
$$
Thus $u\in L^\infty(0,T;W^{1,3})$ and $u\cdot\nabla u\in L^\infty(0,T;L_p^2)$. Repeating the above argument yields
\bg\notag%\label{uh2}
u\in W^{1, q}(W_p^{2,2}), \quad \forall q\in (1, \infty),
\ed
from which it follows that $u\in C^\alpha([0,T]\times Q)$.

Back to equation (\ref{eqL3}),  we conclude that
$$
d \in C^{1+\frac \alpha 2,2+\alpha}((0, T)\times Q), \, \, \mbox { for some } \alpha\in(0,1).
$$
Finally $\rho\in C^1((0, T)\times Q)$ follows from the regularity of $u$ and the regularity of the pressure $p$ follows easily from the regularity of $(\rho, u, d)$ similarly as in \cite{AKM}. This completes the proof of Theorem \ref{Mthm2}.\\

%Proof of Corollary \ref{Cor2}: Follows by Krylov's Theorem \cite{Kr2} just as in the case of two dimension.

%88888888888888888888888888888888888888888888888888888888888888888888888
\bigskip

\section{Appendix: Proof of Uniqueness}
\label{sec:uni}
In this appendix we sketch the proof of Theorem \ref{unique}. For  the full Ericksen-Leslie system with constant density, Wu, Xu and Liu \cite{WXLn} proved that the regular solution $(u,d)$ is unique in the sense: two regular solutions starting from the same initial data are identical. In present paper, for the full Ericksen-Leslie system with non-constant density, we show a weak-strong type of uniqueness result with certain conditions on the weak solution. The idea is to calculate the energy law satisfied by the difference of the regular and weak solutions and establish a Gronwall's inequality. In our case, to calculate the energy law of the difference of the regular and the weak solutions it has some extra terms involving the density. We shall proceed an analogous computation as we did in \cite{DQSr} to achieve the goal.  Similarly, the estimates are more involved  requiring  additional  bounds on the strong solution $(\rho,u,d)$ to yield a Gronwall's inequality. In 2D, we need
\bg\label{more2}
\nabla\rho,\nabla u\in L^\infty((0,T)\times\Omega), \ \ u_t, u\cdot\nabla u\in L^\infty(0,T;L^q(\Omega)), q>2.
\ed
In 3D, we need
\bg\label{more3}
\nabla\rho,\nabla u\in L^\infty((0,T)\times\Omega), \ \ u_t, u\cdot\nabla u\in L^\infty(0,T;L^3(\Omega)).
\ed
With the assumption on data, $\rho_0\in C^1(\bar\Omega)$,  $u_0\in C^{2+\alpha}(\bar\Omega)$ and $d_0\in C^{2+\alpha}(\bar\Omega)$,  the solution $(\rho,u,d)$ from Theorem \ref{Mthm1} or Theorem  \ref{Mthm2} satisfies (\ref{more2}) or (\ref{more3}), respectively .\\

%In fact, if the solution $(\rho,u,d)$ satisfies $u\in L^\infty(0,T;H^1)$ and $d\in L^\infty(0,T;H^2)$ as the condition in \cite{LL}, we also can prove the uniqueness. The reason is that such solution is proven to be regular solution in the sense of Corollary \ref{Cor1} or Corollary  \ref{Cor2} and hence satisfies the estimates in (\ref{more2}) or (\ref{more3}).\\

\textbf {Proof of Theorem \ref{unique}: }First recall that the regular solution $(\rho, u,d)$ from Theorem \ref{Mthm1} or Theorem \ref{Mthm2} satisfies the energy equality:
\begin{equation}\label{engeq}
\begin{split}
&\int_{\Omega}\frac 1 2\rho|u|^2+\frac 1 2|\nabla d|^2+
\rho F(d)dx\\
&+\int_0^T\int_{Q}\mu_4|\nabla u|^2-\frac{1}{\lambda_1}|\Delta d-f(d)|^2+\mu_1|d^TAd|^2+(\mu_5+\mu_6+\frac{\lambda_2^2}{\lambda_1})|Ad|^2dxdt\\
=&\int_{Q}\frac 1 2\rho_0|u_0|^2+\frac 1 2|\nabla d_0|^2+\rho_0F(d_0)dx.
\end{split}
\end{equation}
The density   $\rho$ is the strong solution of  the transport equation, hence it  satisfies that
$$
\int_Q\rho^2dx=\int_Q\rho_0^2dx.
$$
On the other hand,  $\bar\rho$ is a weak solution of the transport equation and $M_1\leq\bar\rho\leq M_2$.   We have by hypothesis that
\bg\label{barrho}
\int_\Omega\bar\rho^2dx\leq\int_\Omega\rho_0^2dx.
\ed
Thus,
\begin{align}\label{rhodiff}
\frac 1 2\int_\Omega|\rho-\bar\rho|^2dx&=\frac 1 2\int_\Omega\rho^2dx+\frac 1 2\int_\Omega\bar\rho^2dx-\int_\Omega\rho\bar\rho dx\\
&\leq\int_\Omega\rho_0^2dx-\int_\Omega\rho\bar\rho dx.\notag
\end{align}
%Recall that $\rho\in C^1([0,T]\times\bar\Omega)$.  We  construct a sequence $\psi_n\in C^\infty([0,T]\times\bar\Omega)$ that
%converges in $C^1([0,T]\times\bar\Omega)$ to $\rho$. This can be done easily with a a mollifier. That is let $\phi \in C^\infty([0,T]\times\bar\Omega)$ satisfy $\int_0^T\int_{\Omega} \phi dx dt =1$ and $\mbox{supp} \phi \subset [0,T]\times\bar\Omega$. Then $\phi_{n}(x) = n^4\phi(nx)$. Let
%$\psi_n = \phi_n \ast \rho$. Now take $\psi_n$ as test functions integrate over $ [0,t)\times\bar\Omega)$
%for any $t \in (0,T)$, we look at solutions in $[0,t)\times\bar\Omega$.
%Thus, taking  $\psi_n\in C_0^\infty([0,t)\times\Omega)$ as test functions yields
%and $C^\infty_0((0,T)\times\Omega)$ is dense in $C^1((0,T)\times\Omega)$,
%we have, for the weak solution $\bar\rho$
%\begin{align}\label{test}
%&\int_\Omega\psi_n(0)\rho_0dx- \int_\Omega\psi_n(t)\rho(t)dx       \\
%&=-\int_0^t\int_\Omega (\psi_n)_t\bar\rho dx dt-
%\int_0^t\int_\Omega\nabla(\psi_n)\bar\rho\bar u dx dt.\notag
%\end{align}
%By the estimates on $\bar\rho$ and $\bar u$, we can pass to the limit in (\ref{test}) as  $n \to\infty$, using Lebesgue dominated convergence theorem. This yields
Since $\rho\in C^1([0,T]\times\bar\Omega)$,  we can take $\rho$ as a test function. Thus, multiplying
\bg\notag
\bar\rho_t+\bar u\cdot\nabla\bar\rho=0
\ed
by $\rho$
and integrating by parts yields
\begin{align}\label{limit}
\int_\Omega\rho_0^2dx-\int_\Omega\rho\bar\rho dx&=-\int_0^t\int_\Omega\bar\rho\rho_t dsdx-\int_0^t\int_\Omega(\bar u\cdot\nabla\rho)\bar\rho dsdx\\
&=\int_0^t\int_\Omega(u\cdot\nabla\rho)\bar\rho dsdx-\int_0^t\int_\Omega(\bar u\cdot\nabla\rho)\bar\rho dsdx.\notag
\end{align}
Here we used again that $\rho$ is a classical solution of the transport equation. Substituting (\ref{limit}) in (\ref{rhodiff}) gives
\begin{align}\label{rhodiff1}
\frac 1 2\int_\Omega|\rho-\bar\rho|^2dx&\leq\int_0^t\int_\Omega\bar\rho(u-\bar u)\nabla\rho dsdx\\
&=\int_0^t\int_\Omega(\rho-\bar\rho)(u-\bar u)\nabla\rho dsdx.\notag
\end{align}
Next, calculate the following term
\begin{align}\notag
&\frac 1 2\int_\Omega\bar\rho|u-\bar u|^2dx+\frac 1 2\int_\Omega|\nabla d-\nabla\bar d|^2dx\\
&=\frac 1 2\int_\Omega(\bar\rho-\rho)|u|^2dx+\frac 1 2\int_\Omega\rho|u|^2dx+\frac 1 2\int_\Omega\bar\rho|\bar u|^2dx-\int_\Omega\bar\rho u\cdot\bar udx\notag\\
&+\frac 1 2\int_\Omega|\nabla d|^2dx+\frac 1 2\int_\Omega|\nabla\bar d|^2dx-\int_\Omega \nabla d\otimes\nabla\bar ddx\notag,
\end{align}
where $\nabla d\otimes\nabla \bar d$ denotes the $3\times 3$ matrix whose $ij$-th entry is given by $\nabla_i d\cdot\nabla_j\bar d$ for $1\leq i,j\leq 3$.\\
Using energy equality (\ref{engeq}) for the regular  solution $(\rho,u,d)$ and inequality (\ref{weakeng})  for the weak  solution $(\bar{\rho},\bar{ u},\bar d)$ combined with the last equation yields
\begin{equation}\label{uddiff}
\begin{split}
&\frac 1 2\int_Q\bar\rho|u-\bar u|^2dx+\frac 1 2\int_Q|\nabla d-\nabla\bar d|^2dx\\
\leq &-\int_0^t\int_Q\mu_4|\nabla u-\nabla\bar u|^2-\frac{1}{\lambda_1}|\Delta d-\Delta\bar d|^2+\mu_1|d^TAd-\bar d^T\bar A\bar d|^2dxdt\\
&-\int_0^T\int_Q\left(\mu_5+\mu_6+\frac{\lambda_2^2}{\lambda_1}\right)|Ad-\bar A\bar d|^2dxdt+\frac 1 2\int_Q(\bar\rho-\rho)|u|^2dx\\
&-\int_Q\bar\rho u\otimes\bar u-\rho_0|u_0|^2dx
-\int_Q \nabla d\otimes\nabla\bar d-|\nabla d_0|^2 dx\\
&-\int_Q \rho F(d)+\bar\rho F(\bar d)-2\rho_0F(d_0)dx
-2\mu_1\int_0^t\int_Qd^TAd\bar d^T\bar A\bar d\,dxdt\\
&-2\mu_4\int_0^t\int_Q\nabla u\nabla\bar udxdt+\frac{2}{\lambda_1}\int_0^t\int_Q \Delta d\Delta\bar d\, dxdt\\
&-2\left(\mu_5+\mu_6+\frac{\lambda_2^2}{\lambda_1}\right)\int_0^t\int_QAd\bar A\bar d\, dxdt\\
&+\frac{1}{\lambda_1}\int_0^T\int_Q|f(d)|^2+|f(\bar d)|^2dxdt
-\frac{2}{\lambda_1}\int_0^t\int_Q\Delta df(d)+\Delta \bar df(\bar d)dxdt.
\end{split}
\end{equation}
%Approximating by  smooth sequences the second and third components of the regular  solutions  to (\ref{LCD}) , and using the approximations as multipliers just as we did for  the first component $\rho$  yields
Since $u,d\in C^{1+\alpha/2,2+\alpha}([0,T]\times\bar Q)$, we can take $u,d$ as test functions for the weak
solution $\bar u,\bar d$. Thus it follows that
\begin{equation}\label{uweak}
\begin{split}
\int_Q\bar\rho u\otimes\bar u-\rho_0|u_0|^2dx=&\int_0^t\int_Q\bar\rho\bar uu_tdxdt+\int_0^t\int_Q\bar\rho\bar u(\bar u\cdot\nabla u)dxdt\\
&+\int_0^t\int_Q\nabla \bar d\nabla \bar d\nabla udxdt+\int_0^t\int_Q\nabla\cdot\bar\sigma udxdt,
\end{split}
\end{equation}
\begin{equation}\label{dweak}
\begin{split}
&\int_Q \nabla d\otimes\nabla\bar d-|\nabla d_0|^2 dx\\
=&-\int_0^T\int_Q\Delta\bar dd_t\, dxdt+\int_0^T\int_Q\Delta d(\bar u\cdot\nabla\bar d)dxdt-\int_0^T\int_Q\Delta d\bar\Omega\bar d\, dxdt\\
&+\frac{\lambda_2}{\lambda_1}\int_0^T\int_Q\Delta d\bar A\bar d\, dxdt+\frac{1}{\lambda_1}\int_0^T\int_Q \Delta d(\Delta \bar d-f(\bar d))\, dxdt.
\end{split}
\end{equation}

Substituting (\ref{uweak}) and (\ref{dweak}) in (\ref{uddiff}) and adding (\ref{rhodiff1}) yields
\begin{equation}\label{engddiff}
\begin{split}
&\frac 1 2\int_Q|\rho-\bar \rho|^2dx+\frac 1 2\int_Q\bar\rho|u-\bar u|^2dx+\frac 1 2\int_Q|\nabla d-\nabla\bar d|^2dx\\
\leq &-\int_0^t\int_Q\mu_4|\nabla u-\nabla\bar u|^2-\frac{1}{\lambda_1}|\Delta d-\Delta\bar d|^2+\mu_1|d^TAd-\bar d^T\bar A\bar d|^2dxdt\\
&-\int_0^T\int_Q\left(\mu_5+\mu_6+\frac{\lambda_2^2}{\lambda_1}\right)|Ad-\bar A\bar d|^2dxdt+\frac 1 2\int_Q(\bar\rho-\rho)|u|^2dx\\
&-\int_Q \rho F(d)+\bar\rho F(\bar d)-2\rho_0F(d_0)dx
+\int_0^t\int_Q(\bar \rho-\rho)(u-\bar u)\nabla\rho \, dxdt\\
&-\int_0^t\int_Q\bar\rho\bar uu_t\, dxdt-\int_0^t\int_Q\bar\rho\bar u\bar u\nabla u\, dxdt-\int_0^t\int_Q\nabla \bar d\nabla\bar d\nabla u\, dxdt\\
&-\int_0^t\int_Qu\nabla\cdot\bar\sigma' \, dxdt+\mu_4\int_0^t\int_Q\nabla\bar u\nabla u\, dxdt+
\int_0^t\int_Q\Delta\bar dd_t\, dxdt\\
&-\int_0^t\int_Q\bar u\nabla\bar d\Delta d\, dxdt+\int_0^t\int_Q\bar \Omega\bar d\Delta d\, dxdt-\frac{\lambda_2}{\lambda_1}\int_0^t\int_Q\bar A\bar d\Delta d\, dxdt\\
&-2\mu_1\int_0^t\int_Qd^TAd\bar d^T\bar A\bar d\,dxdt
-2\mu_4\int_0^t\int_Q\nabla u\nabla\bar udxdt\\
&+\frac{2}{\lambda_1}\int_0^t\int_Q \Delta d\Delta\bar d\, dxdt
-2\left(\mu_5+\mu_6+\frac{\lambda_2^2}{\lambda_1}\right)\int_0^t\int_QAd\bar A\bar d\, dxdt\\
&+\frac{1}{\lambda_1}\int_0^T\int_Q|f(d)|^2+|f(\bar d)|^2dxdt
-\frac{2}{\lambda_1}\int_0^t\int_Q\Delta df(d)+\Delta \bar df(\bar d)dxdt.
\end{split}
\end{equation}
Note that
\begin{equation}\notag
\begin{split}
-\int_0^t\int_Q\bar\rho\bar uu_t\, dxdt
=&-\int_0^t\int_Q(\bar\rho-\rho)(\bar u-u)u_t\, dxdt\\
&-\int_0^t\int_Q(\bar\rho-\rho) uu_t\, dxdt-\int_0^t\int_Q\rho\bar uu_t\, dxdt,
\end{split}
\end{equation}
while 
\begin{equation}\notag
\begin{split}
-\int_0^t\int_Q(\bar\rho-\rho) uu_t\, dxdt=&-\frac{d}{dt}\int_0^t\int_Q(\bar\rho-\rho)\frac{|u|^2}{2}\, dxdt
+\int_0^t\int_Q(\bar\rho-\rho)_t\frac{|u|^2}{2} \, dxdt\\
=&-\int_0^t\int_Q(\bar\rho-\rho)\frac{|u|^2}{2}\, dxdt
-\int_0^t\int_Q\div(\bar\rho\bar u-\rho u)\frac{|u|^2}{2} \, dxdt\\
=&-\int_0^t\int_Q(\bar\rho-\rho)\frac{|u|^2}{2}\, dxdt
+\int_0^t\int_Q\bar\rho\bar uu\nabla u \, dxdt\\
&-\int_0^t\int_Q(\rho-\bar\rho) uu\nabla u \, dxdt-\int_0^t\int_Q\bar\rho uu\nabla u \, dxdt
\end{split}
\end{equation}
and 
\begin{equation}\notag
\begin{split}
-\int_0^t\int_Q\rho\bar uu_t \, dxdt=&\mu_4\int_0^t\int_Q\nabla u\nabla \bar u \, dxdt+
\int_0^t\int_Q\rho u\bar u\nabla u \, dxdt\\
&+\int_0^t\int_Q\nabla d\Delta d\bar u \, dxdt-\int_0^t\int_Q\bar u\nabla\cdot\sigma' \, dxdt.
\end{split}
\end{equation}
On the other hand, we have
\begin{equation}\notag
\begin{split}
&\int_0^t\int_Q\Delta\bar dd_t\, dxdt\\
=&-\int_0^t\int_Q\Delta\bar du\nabla d\, dxdt+\int_0^t\int_Q\Delta\bar d\Omega d\, dxdt-\frac{\lambda_2}{\lambda_1}\int_0^t\int_Q\Delta\bar d Ad\, dxdt\\
&-\frac{1}{\lambda_1}\int_0^t\int_Q\Delta\bar d\Delta d\, dxdt+\frac{1}{\lambda_1}\int_0^t\int_Q\Delta\bar d f(d)\, dxdt.
\end{split}
\end{equation}
Thus with some cancelation (\ref{engddiff}) becomes 
\begin{equation}\label{engddiff1}
\begin{split}
&\frac 1 2\int_Q|\rho-\bar \rho|^2dx+\frac 1 2\int_Q\bar\rho|u-\bar u|^2dx+\frac 1 2\int_Q|\nabla d-\nabla\bar d|^2dx\\
\leq &-\int_0^t\int_Q\mu_4|\nabla u-\nabla\bar u|^2-\frac{1}{\lambda_1}|\Delta d-\Delta\bar d|^2+\mu_1|d^TAd-\bar d^T\bar A\bar d|^2dxdt\\
&-\int_0^T\int_Q\left(\mu_5+\mu_6+\frac{\lambda_2^2}{\lambda_1}\right)|Ad-\bar A\bar d|^2dxdt\\
&+\int_0^t\int_Q(\bar \rho-\rho)(u-\bar u)\nabla\rho \, dxdt-\int_0^t\int_Q(\bar \rho-\rho)(\bar u- u)u_t \, dxdt\\
&-\int_0^t\int_Q\bar\rho\bar u\bar u\nabla u-\bar\rho\bar u u\nabla u+\bar\rho u u\nabla u-\bar\rho\bar u u\nabla u\, dxdt\\
&-\int_0^t\int_Q(\rho-\bar\rho) u u\nabla u-(\rho-\bar\rho) \bar u u\nabla u\, dxdt\\
&+\int_0^t\int_Q\nabla \bar d\Delta\bar d u+\nabla  d\Delta d\bar u-\nabla \bar d\Delta d\bar u-
\nabla  d\Delta\bar d u\, dxdt\\
&-\int_0^t\int_Qu\nabla\cdot\bar\sigma' +\bar u\nabla\cdot\sigma' \, dxdt\\
&-\int_0^t\int_Q2\mu_1d^TAd\bar d^T\bar A\bar d
+2\left(\mu_5+\mu_6+\frac{\lambda_2^2}{\lambda_1}\right)Ad\bar A\bar d\, dxdt\\
&+\int_0^t\int_Q \Omega d\Delta\bar d+\bar \Omega\bar d\Delta d-\frac{\lambda_2}{\lambda_1}(\bar A\bar d\Delta d+Ad\Delta\bar d)\, dxdt\\
&+\frac{1}{\lambda_1}\int_0^t\int_Q\Delta\bar df(d)+\Delta df(\bar d)-2\Delta df(d)-2\Delta \bar df(\bar d)dxdt\\
&-\int_Q \rho F(d)+\bar\rho F(\bar d)-2\rho_0F(d_0)dx+\frac{1}{\lambda_1}\int_0^T\int_Q|f(d)|^2+|f(\bar d)|^2dxdt\\
:=&-G_1-G_2+I_1+I_2+\ldots+I_{11}.
\end{split}
\end{equation}
Note that 
\begin{equation}\notag
\begin{split}
I_3=&-\int_0^t\int_Q\bar\rho \nabla u|u-\bar u|^2dxdt,\\
I_4=&\int_0^t\int_Q(\rho-\bar\rho)u \nabla u(\bar u-u)dxdt,\\
I_5=&\int_0^t\int_Qu(\nabla d-\nabla\bar d)(\Delta d-\Delta\bar d)dxdt\\
&-\int_0^t\int_Q(\bar u-u)(\nabla\bar d-\nabla d)\Delta d\, dxdt.
%I_9=&-\frac{1}{\lambda_1}\int_0^t\int_Q(\Delta d-\Delta\bar d)(f(d)-f(\bar d))dxdt\\
%&-\frac{1}{\lambda_1}\int_0^t\int_Q\Delta df(d)-\Delta\bar df(\bar d)dxdt:=I_{91}+I_{92}.
\end{split}
\end{equation}
Using the relation $\lambda_2=\mu_2+\mu_3$, $I_6+I_7$ can be reorganized as
\begin{equation}\notag
\begin{split}
I_6+I_7=&\mu_1\int_0^t\int_Qd^TAd(d\otimes d-\bar d\otimes\bar d)(\nabla \bar u-\nabla u)dx dt\\
&+\mu_1\int_0^t\int_Q(d^TAd-\bar d^T\bar A\bar d)(d\otimes d-\bar d\otimes\bar d)\nabla udx dt\\
&+\left(\mu_5+\mu_6+\frac{\lambda_2^2}{\lambda_1}\right)\int_0^t\int_Qd^TA(d-\bar d)(\nabla \bar u-\nabla u)dx dt\\
&+\left(\mu_5+\mu_6+\frac{\lambda_2^2}{\lambda_1}\right)\int_0^t\int_Q(d^TA-\bar d^T\bar A)(d-\bar d)\nabla udx dt
\end{split}
\end{equation}
Applying (\ref{eqL3}) and the fact $\int_Qu\cdot\nabla df(d)dx=0$, we derive that
\begin{equation}\label{eq:I8}
\begin{split}
&I_8+I_9+I_{10}+I_{11}\\
=&-\int_0^t\int_Q\left(\Omega-\frac{\lambda_2}{\lambda_1}A\right)(d-\bar d)(\Delta d-\Delta\bar d)dxdt\\
&+\int_0^t\int_Q\left(\left(\Omega-\frac{\lambda_2}{\lambda_1}A\right)-\left(\bar\Omega-\frac{\lambda_2}{\lambda_1}\bar A\right)\right)(d-\bar d)\Delta d\, dxdt\\
&-\frac{1}{\lambda_1}\int_0^t\int_Q(\Delta d-\Delta\bar d)(f(d)-f(\bar d))dxdt
\end{split}
\end{equation}

Recall that the regular solution $(\rho,u,d)$ satisfies (\ref{more2}) in 2D and (\ref{more3}) in 3D.  By H\"older  and Gagliardo-Nirenberg inequalities on the terms of $I_1, \ldots, I_{11}$, it follows that
\begin{align}\label{Grnw}
&\frac 1 2\int_Q|\rho(t)-\bar\rho(t)|^2+\bar\rho(t)|u(t)-\bar u(t)|^2+|\nabla d(t)-\nabla\bar d(t)|^2dx\\
&\leq C\int_0^t\int_Q|\rho-\bar\rho|^2+\bar \rho |u-\bar u|^2+|\nabla d-\nabla\bar d|^2dxdt.\notag
\end{align}
To handle the last integral in (\ref{eq:I8}), we used the fact that $|d|\leq 1$ and the hypothesis $|\bar d|$ is bounded, which imply $|f(d)-f(\bar d)|\leq C|d-\bar d|$ by the definition of $f(d)$. Thus,
\bg\notag
\int_0^t\int_Q|f(d)-f(\bar d)|^2dxdt\leq C\int_0^t\int_Q|d-\bar d|^2dxdt\leq C(Q)\int_0^t\int_Q|\nabla d-\nabla\bar d|^2dxdt
\ed
where the constant $C$ depends on space domain $Q$ not on time $T$.% and $C$ depends on the dimension of the space.\\

Thanks to the fact $\bar\rho\geq M_1>0$, applying Gronwall's inequality to (\ref{Grnw}) we obtain
\begin{align}\notag
&\frac 1 2\int_\Omega|\rho(t)-\bar\rho(t)|^2+\bar\rho|u(t)-\bar u(t)|^2+|\nabla d(t)-\nabla\bar d(t)|^2dx\\
&\leq \int_\Omega|\rho(0)-\bar\rho(0)|^2+\bar\rho(0)|u(0)-\bar u(0)|^2+|\nabla d(0)-\nabla\bar d(0)|^2dxe^{Ct}\notag\\
&=0\notag
\end{align}
for all $t>0$ which implies
$$
\bar\rho-\rho=\bar u-u=\bar d-d\equiv 0.
$$
This completes the proof of Theorem \ref{unique}.

%88888888888888888888888888888888888888888888888888888888888888888888888
%88888888888888888888888888888888888888
\bigskip

%888888888888888888888888888888888888888888888888888888888888888888888

{}

\end{document}